\documentclass[12pt]{article}
\usepackage{amsmath}
\usepackage{amsfonts}
\usepackage{amssymb}
\usepackage{pb-diagram,lamsarrow,pb-lams}
\def\qedsymbol{~QED.}
\def\proofname{Proof.}
\newenvironment{Proof}{\par\noindent{\it\proofname}}{{\unskip\nobreak\hfill{\it\qedsymbol}}\par\vskip 9pt}
\newenvironment{Proof*}{\par\noindent}{{\unskip\nobreak\hfill{\it\qedsymbol}}\par\vskip 9pt}
\ifx\undefined\bysame \newcommand{\bysame}{\leavevmode\hbox to3em{\hrulefill}\,}\fi
\numberwithin{equation}{section}
\newtheorem{Thm}{Theorem}[section]
\newtheorem{Lem}[Thm]{Lemma}
\newtheorem{Def}[Thm]{Definition}
\newtheorem{Prop}[Thm]{Proposition}
\newtheorem{Fact}[Thm]{Fact}
\newtheorem{Cor}{Corollary}[Thm]
\newtheorem{Conj}[Thm]{Conjecture}
\newtheorem{Exam}[Thm]{Example}

\newtheorem{Rem}[Thm]{Remark}
\newtheorem{Rems}[Thm]{Remark}

\newcommand{\cat}{\operatorname{cat}}

\newcommand{\catr}{\operatorname{cat_{\it 3}}}
\newcommand{\catp}{\operatorname{cat_{\it p}}}
\newcommand{\catq}{\operatorname{cat_{\it q}}}

\newcommand{\G}[1]{{\Omega}#1}

\newcommand{\integer}{\mathbb Z}

\newcommand{\real}{\mathbb R}
\newcommand{\complex}{\mathbb C}
\newcommand{\quaternion}{\mathbb H}
\newcommand{\cayley}{{\mathbb O}}

\newcommand{\homeo}{\approx}

\newcommand{\sq}{{\mathcal S\kern-0.2em q}}

\newcommand{\Ext}{\operatorname{Ext}}
\newcommand{\ad}{\operatorname{ad}}
\newcommand{\ads}{\mbox{\scriptsize\rm ad}\,}
\newcommand{\proj}{\operatorname{pr}}
\newcommand{\incl}{\operatorname{in}}
\newcommand{\comp}{\smash{\lower-.1ex\hbox{\scriptsize$\circ$}}}
\title{$A_\infty$-method in Lusternik-Schnirelmann category}
\author{Norio Iwase\thanks{%
Dedicated to Professor Masayoshi Kamata for his 60th birthday.
This research was partially supported by Max-Planck-Institute f\"ur Mathematik and Grant-in-Aid for Scientific Research (C)08640125 from The Ministry of Science, Sports and Culture.%
\endgraf
{\it MSC:\/} Primary 55M30, Secondary 55P35, 55Q25, 55R35, 55S36.
\endgraf
{\it Keywords and phrases.\/} LS category, higher homotopy associativity, homology decomposition, sphere bundles over spheres, manifold counter example to the Ganea conjecture.
}\\
\\
\small {\it Address\/}:
Graduate School of Mathematics,
Kyushu University, Japan.
\\
\small {\it e-mail\/}: 
iwase@math.kyushu-u.ac.jp
}
\date{\small\today}
\begin{document}
\maketitle
%
%
\begin{abstract}\noindent%
To clarify the method behind \cite{Iwase:counter-ls}, a generalisation of Berstein-Hilton Hopf invariants is defined as `higher Hopf invariants'.
They detect the higher homotopy associativity of Hopf spaces and are studied as obstructions not to increase the LS category by one by attaching a cone.
Under a condition between dimension and LS category, a criterion for Ganea's conjecture on LS category is obtained using the generalised higher Hopf invariants, which yields the main result of \cite{Iwase:counter-ls} for all the cases except the case when $p=2$.
As an application, conditions in terms of homotopy invariants of the characteristic maps are given to determine the LS category of sphere-bundles-over-spheres.
Consequently, a closed manifold $M$ is found not to satisfy Ganea's conjecture on LS category and another closed manifold $N$ is found to have the same LS category as its `punctured submanifold' $N-\{P\}, P \in N$.
But all examples obtained here support the conjecture in \cite{Iwase:counter-ls}.
\end{abstract}
%
%
\baselineskip 21pt
\section{Introduction}
\par\noindent
In this paper, each space is assumed to have the homotopy type of a CW complex.
The LS category of $X$ is the least number $m$ such that there is a covering of $X$ by $m+1$ open subsets each of which is contractible in $X$, which is (by Whitehead \cite{Whitehead:elements}) the least number $m$ such that the diagonal map $\Delta_{m+1} : X \to X^{m+1}$ can be compressed into the `fat wedge' $T^{m+1}(X)$ or $X^{[m+1]}$.
Hence $\cat{\{\ast\}} = 0$.

As is well-known, the LS category of a product space $\cat{X{\times}S^n}$ is either $\cat{X}$ or $\cat{X}+1$.
A problem was posed by Ganea in \cite{Ganea:conjecture}: Can only the latter case occur for any $X$ and $n \geq 1$?

The affirmative answer had been supposed to be true and came to be known as `Ganea's conjecture' (see \cite{Hess:ls-cat_0}) or `the Ganea conjecture' (see \cite{James:ls-category-survey}).
A major advance in this subject was made by Jessup \cite{Jessup:ls-cat_0} and Hess \cite{Hess:ls-cat_0} working in the rational category: the rational version of the conjecture is true for $n \geq 2$.
Also by Singhof \cite{Singhof:minimal-cover} and Rudyak \cite{Rudyak:ls-cat_mfds}, \cite{Rudyak:ls-cat_mfds2}, the conjecture is true for a large class of manifolds.

%
However in June 1997, the author found a counter example (see \cite{Iwase:counter-ls_o}), in an effort to provide a criterion for establishing the conjecture (which is given in this paper as Theorems \ref{thm:stably-trivial}, \ref{thm:stably-non-trivial} and Corollary \ref{cor:criterion-ls}), using properties of higher Hopf invariants (see \cite{Iwase:counter-ls-m_o}) and fibrations associated with the $A_\infty$-structure of $\G{X}$ (see Sugawara \cite{Sugawara:hopf}, \cite{Sugawara:group-like}, Stasheff \cite{Stasheff:higher-associativity} and Iwase-Mimura \cite{IM:higher-associativity}).
The author knows that Don Stanley was trying to find out a counter example using the ordinary James-Hopf invariants, and also Lucile Vandembroucq \cite{Vandembroucq:sigma-category} obtained a related result on a sufficient condition to Ganea's conjecture at about the same time.
The author also knows that soon after \cite{Iwase:counter-ls-m_o}, the higher Hopf invariants were begun to be studied by Stanley (see \cite{Stanley:cat-cl}).

This paper is organised as follows:
In Section \ref{sect:proj-hopf}, to clarify the method behind \cite{Iwase:counter-ls}, a generalisation of the Berstein-Hilton Hopf invariants is defined with its related invariants to detect the higher homotopy associativity of a Hopf space.
In Section \ref{sect:hopf-cat}, under a condition between dimension and LS category, some conjectures on LS category are verified by using fibrations associated with the $A_{\infty}$-structure of a loop space.
In Section \ref{sect:hopf-diag}, a result of Boardman-Steer is generalised to give a sufficient condition to determine LS category in terms of a generalised version of the Berstein-Hilton crude Hopf invariants.
In Section \ref{sect:hom-decomp}, the relation between a homology decomposition and LS category of (product) spaces is shown, by extending a result of Curjel \cite{Curjel:k'inv_co-hopf}.
In Section \ref{sect:hhi-example}, generalising the main result of \cite{Iwase:counter-ls} for all the cases except the case when $p=2$, some more examples are obtained by the properties of the higher Hopf invariants given in Section \ref{sect:proj-hopf}.
In Section \ref{sect:cat-sphere-bundle-over-sphere}, we give some conditions to determine the LS category of sphere-bundles-over-spheres.
Using it, we construct, in Section \ref{sect:mfd-example}, an orientable closed manifold $N_p$, for each prime $p \geq 5$, with the LS category same as its `punctured submanifold' $N_p - \{P\}$, $P \in N_p$.
Also another orientable closed manifold $M$ is constructed as a counter example to Ganea's conjecture.

The author would like to express his gratitude to Ioan James for giving an attractive lecture on LS category at the University of Aberdeen which inspired him to consider LS category again, John Hubbuck, Yuly Rudyak, Kouyemon Iriye, Donald Stanley, Daniel Tanr{\'e}, Octavian Cornea, Hans Scheerer and Hans Baues for valuable conversations and encouragement which helped to organise his thoughts, the University of Aberdeen for its hospitality during the author's stay in 1997, Max-Planck-Institut f\"ur Mathematik for its hospitality during the author's stay in 2000 and the members of the Graduate School of Mathematics Kyushu University for allowing him to be away for a long term, without which this work could not be done.

\section{Projective spaces and higher Hopf invariants}\label{sect:proj-hopf}

In this section, we introduce a generalised version of the Berstein-Hilton Hopf invariant (see \cite{BH:category}), a higher Hopf invariant for short, in terms of projective spaces associated with the $A_{\infty}$-structure of a loop space, to detect the higher homotopy associativity, or the $A_m$-structure of a Hopf space (see Example \ref{exam:A_m-invariant-one} and Conjecture \ref{conj:A_m-invariant-one}).
We also show that a higher Hopf invariant gives the obstruction for increasing the LS category by one by attaching a cone, as in \cite{Iwase:counter-ls}.

For a given space, its loop space is an $A_{\infty}$-space with the given space as its $A_{\infty}$-structure.
More precisely, every space $X$ has a filtration given by the projective spaces $P^m(\Omega{X})$ of its loop space $\Omega{X}$.
There is a ladder of Stasheff's fibrations $E^{m+1}(\Omega{X}) \overset{p^{\G{X}}_{m}}\to P^m(\Omega{X})$ with the fibre $\Omega{X}$ contractible in $E^{m+1}(\Omega{X})$, if $m \geq 1$.
The total space $E^{m+1}(\Omega{X})$ has the homotopy type of the $m+1$-fold unreduced join of $\Omega{X}$ which is denoted by $\bar{E}^{m+1}(\Omega{X})$ and the base space $P^{m+1}(\Omega{X})$ has the homotopy type of the mapping cone of $p^{\G{X}}_{m}$, $m \geq 0$: the fibration is induced by the inclusion $e^X_m : P^m(\Omega{X}) \to P^{\infty}(\Omega{X})$ from the {\it universal} fibration $E^{\infty}(\Omega{X}) \overset{p^{\G{X}}_{\infty}}\to P^{\infty}(\Omega{X})$ with contractible total space. (see \cite{Stasheff:higher-associativity} for details.)

\begin{Thm}[Ganea]\label{thm:domination}
Let $X$ be a connected CW complex.
Then $\cat{X} \leq m$ if and only if the inclusion $e^X_m : P^{m}(\Omega{X}) \subset P^{\infty}(\Omega{X}) \simeq X$ has a right homotopy inverse (homotopy section).
\end{Thm}
This result enables us to define local versions of $\cat{}$, e.g, $\catp$ is defined in \cite{Iwase:counter-ls}, for a prime $p$, as the least number $m$ such that $e^X_m : P^{m}(\Omega{X}) \subset P^{\infty}(\Omega{X}) \simeq X$ has a homotopy section at $p$, e.g. $\catp{S_{(p)}^1} = 2$ while $\cat{S^1} = 1$.
We remark that in some of the literature, the composition functor $P^m{\comp}\Omega$ is abbreviated as $G_m$, the `Ganea space' functor.
We use the following fact.
\begin{Fact}[\cite{Iwase:counter-ls}]\label{cor:product-formula-sphere}
Let $X$ be a connected CW complex with $\cat{X} = m$.
Then $\cat{X{\times}S^n} = m$ if and only if $(e^X_m{\times}1_{S^n})\vert : P^{m}(\Omega{X}){\times}\{\ast\} \cup P^{m-1}(\Omega{X}){\times}S^n \subset P^{\infty}(\Omega{X}){\times}S^n \simeq X{\times}S^n$ has a homotopy section.
\end{Fact}
In \cite{Iwase:counter-ls}, a conjecture was posed instead of Ganea's conjecture (Conjecture 1.4 in \cite{Iwase:counter-ls}):
\begin{Conj}[\cite{Iwase:counter-ls}]\label{conj:iwase-number}
For any space $X$, there exists an integer $n(X)$, $1 \leq n(X) \leq \infty$, such that $\cat{X{\times}S^n}$ is equal to $\cat{X}+1$, if $n < n(X)$; $\cat{X}$, if $n \geq n(X)$.
\end{Conj}

When $V$ is the suspension of a co-H-space, say $V = \Sigma{V_0}$ with $V_0$ a co-H-space, we fix a canonical structure map $\sigma(V) = \Sigma\ad(1_{V}) : V \to \Sigma\G{V} = P^1(\G{V})$ for $\cat{V} \leq 1$, i.e. a (homotopy) section of the evaluation map $e^V_1$, where $\ad(1_{V}) : V_0 \to \G{\Sigma{V_0}}$ is the adjoint of the identity.
Then $\sigma(V)$ gives a homotopy commutative and homotopy associative co-H-structure on $V$.
In this section, we fix a non-contractible co-H-space $V$ with right and left inversion (e.g, $V$ is a suspension space), together with a structure map $\sigma(V) : V \to {\Sigma}{\Omega}V$ for $\cat{V} = 1$.
For any given $m \geq 1$, we often regard $P^m(\G{V})$ as the target of $\sigma(V)$, since ${\Sigma}{\Omega}V = P^1(\G{V}) \subset P^m(\G{V})$.

\begin{Def}\label{def:higher-hopf}\quad
Let $X$ be a space with $\cat{X} \leq m$, $m \geq 1$.
For a choice of the homotopy section $\sigma(X) : X \to P^{m}({\Omega}X)$ ($m \geq 1$), we define a higher Hopf invariant as
\begin{equation*}
H_m = H^{\sigma(X)}_m : [V,X] \to [V,E^{m+1}(\Omega{X})],
\end{equation*}
which is a homomorphism when $V$ is homotopy associative and homotopy commutative:
For a map $f : V \to X$, the difference between $\sigma(X){\comp}f$ and $P^m(\G{f}){\comp}\sigma(V)$, in the algebraic loop $[V,P^{m}(\Omega{X})]$, is given by a map $d^{\sigma(X)}_m(f) : V \to P^{m}(\Omega{X})$ so that $\sigma(X){\comp}f + d^{\sigma(X)}_m(f) \sim P^m(\G{f}){\comp}\sigma(V)$.
Since $d^{\sigma(X)}_m(f)$ vanishes in $[V,P^{\infty}(\Omega{X})] \cong [V,X]$, it has a unique lift $H^{\sigma(X)}_m(f) : V \to E^{m+1}(\Omega{X})$ to the total space of Stasheff's fibration
\begin{equation}\label{eq:Stasheff-fibration}
\G{X} \hookrightarrow E^{m+1}(\G{X}) \overset{p^{\G{X}}_{m}}\to P^{m}(\G{X}) \overset{e^X_{m}}\to X, \quad m \geq 1.
\end{equation}
\end{Def}
\begin{Rems}\label{rem:higher-hopf}
\begin{enumerate}
\item\label{rems:(1)}\quad
When $V$ is a homotopy associative co-H-space, $\sigma(V)$ is a co-H-map by Theorem 2.2 of Ganea \cite{Ganea:co-h-space}.
Hence a simple calculation shows that $d^{\sigma(X)}_m$ and $H^{\sigma(X)}_m$ are homomorphisms, if $V$ is a homotopy associative and homotopy commutative co-H-space.
\item\label{rems:(2)}\quad
By Berstein-Dror \cite{BD:ass-co-act}, a homotopy associative co-H-space admits another homotopy associative co-H-structure which has right and left inversions.\par
\item\label{rems:(3)}\quad
When $V$ is a Moore space of type $(A,n)$, $H_m$ can be regarded as the Berstein-Hilton Hopf invariant $H : \pi_n(X;A) \to \pi_n(X^{m+1},T^{m+1}X;A)$ (see \cite{BH:category}), since $\pi_n(X^{m+1},T^{m+1}X;A) = \pi_n(E^{m+1}(\G{X});A)$ by Ganea (see the proof of Theorem 1.1 in \cite{Iwase:counter-ls}).\par
\item\label{rems:(4)}\quad
When $V = \Sigma{V_0}$ is a suspension space, a map $f : V \to X$ factors through $\Sigma\G{X}$ as $f = e^X_1{\comp}{\Sigma}\ad(f)$, where $e^X_1$ denotes the evaluation map and $\ad(f) = \G{f}{\comp}\ad(1_V) : V_0 \to \G{X}$ is the adjoint of $f$.
Then the higher Hopf invariant $H^{\sigma(V)}_m(e^X_1)$ of $e^X_1 : \Sigma\G{X} \to X$ has a kind of `universality': $H^{\sigma(V)}_m(f) = H^{\sigma(\Sigma\G{X})}_m(e^X_1){\comp}{\Sigma}\ad(f)$ (see Corollary \ref{cor:universal_hopf}).
So we call $H^{\sigma(\Sigma\G{X})}_m(e^X_1)$ the \underline{universal Hopf invariant} for $X$ and $m \geq \cat{X}$.
\par
\item\label{rems:(5)}\quad
If $X$ is a suspension space, say $X = \Sigma{Y}$, then we have $P^1(\G{X}) = \Sigma{\G{X}} = \Sigma{J(Y)} \simeq \Sigma{Y} \vee \Sigma{\overset{2}{\wedge}Y} \vee \Sigma{\overset{3}{\wedge}Y} \vee ...$.
Let us recall that $H^{\sigma(X)}_1$ is uniquely determined by $d^{\sigma(X)}_1$, whose projection to $\Sigma{\overset{j}{\wedge}Y}$ gives exactly the $j$-th James Hopf invariant $h_j$, $j \geq 2$.
Thus we may regard $H_1$ as the collection of all James Hopf invariants $h_j$, $j \geq 2$.\par
\end{enumerate}
\end{Rems}

The above definition of a higher Hopf invariant also allows us to define a generalisation of the Berstein-Hilton crude Hopf invariant as follows.

\begin{Def}\label{def:crude-hopf-invariant}
\begin{equation*}
\dgHORIZPAD=.5em
\dgVERTPAD=1ex
\begin{diagram}
\node{\bar{H}_m = \bar{H}^{\sigma(X)}_m : [V,X]} \arrow{e,t}{H^{\sigma(X)}_m} \node{[V,E^{m+1}(\Omega{X})]} \arrow{e,t}{({\Sigma}{h^X_{m+1}}_{\ast}){\comp}{\Sigma}} \node{[{\Sigma}V,\wedge^{m+1}{\Sigma}\Omega{X}]} \arrow{e,t}{(\wedge^{m+1}e^{X}_{1})_{\ast}} \node{[{\Sigma}V,\wedge^{m+1}{X}],}
\end{diagram}
\end{equation*}
where $e^{X}_{1}$ is the evaluation map and $h^X_{m+1} : E^{m+1}(\Omega{X}) \to \bar{E}^{m+1}(\Omega{X}) \simeq \Omega{X}{\wedge}(\wedge^{m}{\Sigma}\Omega{X})$ denotes the natural homotopy equivalence (see Stasheff \cite{Stasheff:higher-associativity}).
\end{Def}

\begin{Exam}\label{exam:A_m-invariant-one}
For an $A_m$-space $G$ in the sense of Stasheff \cite{Stasheff:higher-associativity}, the adjoint of the inclusion $\iota^G_{1,m} : {\Sigma}G \hookrightarrow P^m(G)$ is an $A_m$-map $\ad(\iota^G_{1,m}) : G \hookrightarrow \G{P^m(G)}$ whose $A_m$-structure map $P^m(\ad(\iota^G_{1,m})) : P^m(G) \to P^m(\G{P^m(G)})$ is given by a splitting of $e^{P^m(G)}_m : P^{m}(\G{P^m(G)}) \to P^m(G)$ (see \cite{IM:higher-associativity}).
By putting $X = P^{m}(G)$, $V = \bar{E}^{m+1}(G)$ and $\sigma(X)=P^m(\ad(\iota^G_{1,m}))$,  we have 
\begin{equation}
H_m = H^{\sigma(X)}_m : [\bar{E}^{m+1}(G),P^{m}(G)] \to [\bar{E}^{m+1}(G),\bar{E}^{m+1}(\G{P^{m}(G)})].
\end{equation}
The group $[\bar{E}^{m+1}(G),\bar{E}^{m+1}(\G{P^{m}(G)})]$ contains the image of $1 = 1_{\bar{E}^{m+1}(G)}$ under the homomorphism $\ad(\iota^G_{1,m})_{\ast} : [\bar{E}^{m+1}(G),\bar{E}^{m+1}(G)] \to [\bar{E}^{m+1}(G),\bar{E}^{m+1}(\G{P^{m}(G)})]$.
As is clearly seen, if $G$ is an $A_{m+1}$-space, then there is a `higher Hopf invariant one' element, i.e. there is a map $f : \bar{E}^{m+1}(G) \to P^{m}(G)$ such that $H_m(f)$ is the image of $1$ as mentioned above.
The converse is also verified when $G$ is a sphere: Let $f : \bar{E}^{m+1}(G) \to P^{m}(G)$ be a `higher Hopf invariant one' element.
Then an easy calculation of the Serre spectral sequence shows that the fibre of $f$ has the homotopy type of the sphere $G$.
Thus  $G$ must be an $A_{m+1}$-space (see \cite{Stasheff:higher-associativity}).
\end{Exam}

This suggests the following conjecture.

\begin{Conj}\label{conj:A_m-invariant-one}
An $A_m$-space $G$ has an $A_{m+1}$-structure extending the given $A_m$-structure if and only if there is a map $f : \bar{E}^{m+1}(G) \to P^{m}(G)$ with `higher Hopf invariant one', where $P^{m}(G)$ is the $G$-projective $m$-space associated with the $A_m$-structure in the sense of Stasheff.
\end{Conj}

We bring this higher Hopf invariant $H^{\sigma(X)}_m$ into stable homotopy theory:
\begin{Def}\label{defn:stabilised-higher-hopf}
We define a stabilised higher Hopf invariant as
\begin{equation*}
{\mathcal H}_m  = {\mathcal H}^{\sigma(X)}_m : [V,X] \overset{H^{\sigma(X)}_m}\to [V,E^{m+1}(\Omega{X})] \overset{{\Sigma^{\infty}_{\ast}}}\to \{V,E^{m+1}(\Omega{X})\}
\end{equation*}
and the stabilised crude higher Hopf invariant as
\begin{equation*}
\bar{\mathcal H}_m  = \bar{\mathcal H}^{\sigma(X)}_m : [V,X] \overset{\bar{H}^{\sigma(X)}_m}\to [{\Sigma}V,\wedge^{m+1}{X}] \overset{{\Sigma^{\infty}_{\ast}}}\to \{{\Sigma}V,\wedge^{m+1}{X}\}.
\end{equation*}
\end{Def}

These definitions of higher Hopf invariants depend on the choice of the structure map $\sigma(X)$.
So it might be useful to define the following set-valued functions.
\begin{Def}\label{def:set-higher-hopf}\quad
\begin{align*}&
H^S_m(f) = \{H^{\sigma_(X)}_{m}(f) \,\vert\, \text{$\sigma(X)$ is a structure map for $\cat{X} = m$}\} \subset [V,E^{m+1}(\G{X})],
\\&
\bar{H}^S_m(f) = \{\bar{H}^{\sigma_(X)}_{m}(f) \,\vert\, \text{$\sigma(X)$ is a structure map for $\cat{X} = m$}\} \subset [{\Sigma}V,\wedge^{m+1}{X}],
\\&
{\mathcal H}^S_m(f) = \{{\mathcal H}^{\sigma_(X)}_{m}(f) \,\vert\, \text{$\sigma(X)$ is a structure map for $\cat{X} = m$}\} \subset \{V,E^{m+1}(\G{X})\}
\\&
\bar{\mathcal H}^S_m(f) = \{\bar{\mathcal H}^{\sigma_(X)}_{m}(f) \,\vert\, \text{$\sigma(X)$ is a structure map for $\cat{X} = m$}\} \subset \{{\Sigma}V,\wedge^{m+1}{X}\}.
\end{align*}
\end{Def}

We show the fundamental properties of higher Hopf invariants.
\begin{Prop}\label{prop:fundamental-property}
Let $V$, $X$ and $f$ be as above.
Then the following two statements hold.
\begin{description}
\item{(1)}
Let $V'$ be the suspension space of a co-H-space.
If $g : V' \to V$ is a co-H-map (or equivalently, $P^1(\G{g}){\comp}\sigma(V') \sim \sigma(V){\comp}g$ in $P^1(\G{V})$), then $H^{\sigma(X)}_m(f{\comp}g) \sim H^{\sigma(X)}_m(f){\comp}g$.
\item{(2)}
Let $X'$ be a space of LS category $\leq m$ with a structure map $\sigma(X')$.
If $h : X \to X'$ is $m$-primitive in the sense of Berstein-Hilton \cite{BH:category} (or equivalently, $\sigma(X'){\comp}h \sim P^m(\G{h}){\comp}\sigma(X)$), then $H^{\sigma(X')}_m(h{\comp}f) \sim E^{m+1}(\G{h}){\comp}H^{\sigma(X)}_m(f)$.
\end{description}
\end{Prop}
\begin{Cor}\label{cor:universal_hopf}
For any map $f : V \to X$, we have the following homotopy relation.
\begin{equation}
H^{\sigma(X)}_m(f) \sim H^{\sigma(X)}_m(e^{X}_1){\comp}\Sigma\G(f){\comp}\sigma(V).
\end{equation}
\end{Cor}
\begin{Rem}
The statement (2) in Proposition \ref{prop:fundamental-property} is a generalisation of Proposition 3.2 in Berstein-Hilton \cite{BH:category}.
\end{Rem}
\begin{Proof*}{\it Proof of Proposition \ref{prop:fundamental-property}.}
Firstly we show (1): Let $d^{\sigma(X)}_m(f)$ be the difference between $\sigma(X){\comp}f$ and $P^m(\G{f}){\comp}\sigma(V)$, in the algebraic loop $[V,P^{m}(\Omega{X})]$.
Since $g$ is a co-H-map, we obtain the following equation up to homotopy:
\begin{equation*}
\sigma(X){\comp}(f{\comp}g) + d^{\sigma(X)}_m(f){\comp}g
\sim \{(\sigma(X){\comp}f) + d^{\sigma(X)}_m(f)\}{\comp}g
\sim P^m(\G{f}){\comp}\sigma(V){\comp}g.
\end{equation*}
Again, since $g$ is a co-H-map, $P^m(\G{g}){\comp}\sigma(V') \sim \sigma(V){\comp}g$, and hence we proceed as
\begin{align*}
\sigma(X){\comp}(f{\comp}g) + d^{\sigma(X)}_m(f){\comp}g
&\sim P^m(\G{f}){\comp}P^m(\G{g}){\comp}\sigma(V') = P^m(\G{(f{\comp}g)}){\comp}\sigma(V').
\end{align*}
This implies the difference between $\sigma(X){\comp}(f{\comp}g)$ and $P^m(\G{(f{\comp}g)}){\comp}\sigma(V')$ is given by $d^{\sigma(X)}_m(f{\comp}g) \sim d^{\sigma(X)}_m(f){\comp}g$, and hence $H^{\sigma(X)}_m(f{\comp}g) \sim H^{\sigma(X)}_m(f){\comp}g$.

Secondly we show (2): Let $d^{\sigma(X)}_m(f)$ be the difference between $\sigma(X){\comp}f$ and $P^m(\G{f}){\comp}\sigma(V)$.
Since $h$ is an $m$-primitive map, it follows that $\sigma(X'){\comp}h \sim P^m(\G{h}){\comp}\sigma(X)$.
Hence we obtain the following equation up to homotopy:
\begin{align*}&
\sigma(X'){\comp}(h{\comp}f) + P^m(\G{h}){\comp}d^{\sigma(X)}_m(f) \sim P^m(\G{h}){\comp}\{\sigma(X){\comp}f + d^{\sigma(X)}_m(f)\}
\\&\qquad
\sim P^m(\G{h}){\comp}P^m(\G{f}){\comp}\sigma(V) = P^m(\G{(h{\comp}f)}){\comp}\sigma(V).
\end{align*}
Thus the difference between $\sigma(X'){\comp}(h{\comp}f)$ and $P^m(\G{(h{\comp}f)}){\comp}\sigma(V)$ is given by $d^{\sigma(X')}_m(h{\comp}f) \sim P^m(\G{h}){\comp}d^{\sigma(X)}_m(f)$, and hence $H^{\sigma(X')}_m(h{\comp}f) \sim P^m(\G{h}){\comp}H^{\sigma(X)}_m(f)$.
\end{Proof*}

\section{Higher Hopf invariant and LS category}\label{sect:hopf-cat}

In the remainder of this paper, we always assume that $m \geq 1$ and $V=\Sigma{V_0}$ a suspension space and we fix the structure map $\sigma(V)=\Sigma\ad(1_V)$, unless otherwise stated.
We begin this section with the following results by James \cite{James:ls-category} and by Berstein and Hilton (see Proposition 2.5 in \cite{BH:category}).

\begin{Prop}[James]\label{prop:dimension-category}
Let a CW complex $X$ be ($d-1$)-connected, $d \geq 2$.
Then the following inequality holds:
$\dim{X} \geq d\cdot\cat{X}$.
\end{Prop}

\begin{Prop}[Berstein-Hilton]\label{prop:uniqueness}
Let a CW complex $X$ be ($d-1$)-connected, $d \geq 2$.
If $\dim{X} \leq {d}m + d - 2$, then a structure map $\sigma(X)$ $:$ $X$ $\to$ $P^m(\G{X})$ for $\cat{X} \leq m$ is determined uniquely up to homotopy.
In particular, when $X$ is simply connected and $\dim{X} \leq {2}\cat{X}$, the structure map $\sigma(X)$ $:$ $X$ $\to$ $P^m(\G{X})$ for $\cat{X} \leq m$ is determined uniquely up to homotopy.
\end{Prop}

Proposition \ref{prop:dimension-category} implies that quite a few complexes satisfy the hypothesis of Proposition \ref{prop:uniqueness}.
Such spaces satisfy the uniqueness of higher Hopf invariants as follows:

\begin{Cor}\label{cor:uniqueness}
Let CW complexes $V$ and $X$ be ($e-1$)-connected and ($d-1$)-connected resp., $d, e \geq 2$ with $\dim{X} \leq {d}m + d - 2$, $\cat{X} \leq m$.
Then the higher Hopf invariant $H_m = H^{\sigma(X)}_m$ is uniquely determined.
\end{Cor}
\begin{Proof}
By Proposition \ref{prop:uniqueness}, the hypothesis on $\dim{X}$ and $\cat{X}$ implies that the structure map $\sigma(X)$ for $\cat{X} \leq m$ is uniquely determined.
By the definition of the higher Hopf invariant $H^{\sigma(X)}_m$, $H_m = H^{\sigma(X)}_m$ is uniquely determined.
\end{Proof}

We describe the relationship between the higher Hopf invariant and the LS category:

\begin{Prop}\label{prop:obstruction}
For a given structure map $\sigma(X)$ for $\cat{X} = m$ and a given map $f : V \to X$, let $W$ be the mapping cone of $f$.
Then the following diagram without the dotted arrows commutes up to homotopy. 
\begin{equation*}
\begin{diagram}
\dgHORIZPAD=.5em
\dgVERTPAD=1ex
\node{V}
        \arrow{e,t}{f}
        \arrow{s,l}{H_{m}(f)}
\node{X}
        \arrow{e,t,J}{i}
        \arrow{s,r}{\sigma(X)}
\node{W}
        \arrow[2]{s,r,..}{\lower-6.5ex\hbox{${}_{\sigma'(W)}$}}
        \arrow{e,t}{q}
\node{{\Sigma}V}
\\
\node{E^{m+1}(\G{X})}
        \arrow{e,b,..}{p^{\G{X}}_{m}}
        \arrow{s,l,L}{E^{m+1}(\G{i})}
\node{P^{m}(\G{X})}
        \arrow{s,r,L}{P^{m}(\G{i})}
\\
\node{E^{m+1}(\G{W})}
        \arrow{e,b}{p^{\G{W}}_{m}}
\node{P^{m}(\G{W})}
        \arrow{e,b,J}{\iota^{\G{W}}_{m}}
\node{P^{m+1}(\G{W})}
        \arrow{e,b}{e^{W}_{m+1}}
\node{W,}
\end{diagram}
\end{equation*}
where 
$i : X \hookrightarrow W$ and $\iota^{\G{W}}_{m} : P^{m}(\G{W}) \to P^{m+1}(\G{W})$ denote the inclusions.
\end{Prop}
\begin{Proof}
Let us recall that $H_{m}(f)$ is given by the unique lift of the difference between $\sigma(X){\comp}f$ and $P^{m}(\G{f}){\comp}\sigma(V)$ in the group $[V,P^m(\G{X})]$.
Hence the composition $P^{m}(\G{i}){\comp}p^{\G{X}}_m{\comp}H_{m}(f)$ gives the difference between $P^{m}(\G{i}){\comp}\sigma(X){\comp}f$ and $P^{m}(\G{i}){\comp}P^{m}(\G{f}){\comp}\sigma(V)$ = $P^{m}(\G{(i{\comp}f)}){\comp}\sigma(V) \sim 0$.
Thus we obtain the commutativity of the diagram.
\end{Proof}
\begin{Rem}\label{rem:category-W'}\quad
By the homotopy commutativity of the left rectangle of the diagram, there is a map $\sigma'(W) : W \to P^{m+1}(\G{W})$ making the diagram commutative, which is given by the homotopy deforming $P^{m}(\G{i}){\comp}\sigma(X){\comp}f$ to $p^{\G{W}}_m{\comp}E^{m+1}(\G{i}){\comp}H_{m}(f)$ in $P^{m}(\G{W})$ and by the map $\chi_{p^{\G{W}}_{m}}{\comp}C(E^{m+1}(\G{i}){\comp}H_{m}(f))$ in $P^{m+1}(\G{W})$, where we denote by $\chi_f : (C(V),V) \to (W,X)$ the relative homeomorphism and by $C$ the functor taking cones, since both top and bottom rows except for the map $e^{W}_{m+1}$ are cofibration sequences.
\end{Rem}
\par
Now we discuss the naturality of $\sigma'(W)$ which is determined as above.
A suspension map $g={\Sigma}g_0 : V' \to V$ between suspension spaces gives the equality
\begin{equation*}
\sigma(V){\comp}g = \Sigma\G{g}{\comp}\sigma(V').
\end{equation*}
For maps $f : V \to X$ and $f' : V' \to X$, we have their (reduced) cofibres $i : X \hookrightarrow W = \{\ast\} \cup [0,1]{\times}V \cup_{f} X$ and $i': X \hookrightarrow W' = \{\ast\} \cup [0,1]{\times}V' \cup_{f'} X$ with null-homotopies $F : 0 \sim i{\comp}f$ and $F' : 0 \sim i'{\comp}f'$:
\begin{equation*}
F(t,v) = t{\wedge}v, \quad F'(t',v') = t'{\wedge}v',\quad\text{for $v \in V$, $v' \in V'$ and $t,t' \in [0,1]$.}
\end{equation*}
The following lemma will be applied in Section \ref{sect:mfd-example}.
\begin{Lem}\label{prop:obstruction-natural}
Let $g : V' \to V$ be a suspension map between suspension spaces.
If $f' \sim f{\comp}g$ with a homotopy $L : [0,1]{\times}V' \to X$, there is a map $\tilde{L} : W' \to W$ extending $i : X \rightarrow W$.
Moreover there is a homotopy $\sigma'(W){\comp}\tilde{L} \sim P^{m+1}(\G{\tilde{L}}){\comp}\sigma'(W')$ relative to $X$, where $\sigma'(W)$ and $\sigma'(W')$ are as in Remark \ref{rem:category-W'} using $H_m(f)$ and $H_m(f')$, resp.
\end{Lem}
\begin{Proof}
Let $V={\Sigma}V_0$ and $V'={\Sigma}V'_0$ and let $g={\Sigma}g_0$, $g_0 : V'_0 \to V_0$.
We define $\tilde{L}$ by the following equation.
\begin{equation*}
\tilde{L}(x) = x,\quad \tilde{L}(t{\wedge}v') = \begin{cases} (2t){\wedge}g(v') \in W,\quad t \leq \frac{1}{2}, \\ L(2-2t,v') \in X \subset W,\quad t \geq \frac{1}{2} \end{cases}
\end{equation*}
for $x \in X$ and $(t,v') \in [0,1]{\times}V'$.
Let us consider the following diagram of homotopies:
\begin{equation}\label{eq:homotopy-of-homotopy}
\divide\dgARROWLENGTH by2
\dgHORIZPAD=.5em
\dgVERTPAD=1ex
\begin{diagram}
\node{P^{m}(\G{(\tilde{L}{\comp}i')}){\comp}\sigma(X){\comp}f'}
	\arrow[2]{e,t,..}{P^{m}(\G{i}){\comp}\sigma(X){\comp}L}
	\arrow[2]{s,..}
\node{}
\node{P^{m}(\G{i}){\comp}\sigma(X){\comp}f{\comp}g}
	\arrow[2]{s,..}
\\
\node{}
\node{(A)}
\node{}
\\
\node{P^{m}(\G{(\tilde{L}{\comp}i')}){\comp}\sigma(X){\comp}f' - 0}
	\arrow[2]{e,t,..}{P^{m}(\G{i}){\comp}\sigma(X){\comp}L - c(0)}
	\arrow[2]{s,l,..}{\substack{c(P^{m}(\G(\tilde{L}{\comp}i')){\comp}\sigma(X){\comp}f') \\ - P^{m}(\G(\tilde{L}{\comp}F')){\comp}\iota_{V'}{\comp}\sigma(V')}}
\node{}
\node{P^{m}(\G{i}){\comp}\sigma(X){\comp}f{\comp}g - 0}
	\arrow[2]{s,r,..}{\substack{c(P^{m}(\G{i})){\comp}\sigma(X){\comp}f{\comp}g \\ - P^{m}(\G{F}){\comp}\iota_{V'}{\comp}\sigma(V')}}
\\
\node{}
\node{(B)}
\node{}
\\
\node{\begin{array}{c}P^{m}(\G{(\tilde{L}{\comp}i')}){\comp}\sigma(X){\comp}f'\\ - P^{m}(\G{(\tilde{L}{\comp}i'{\comp}f')}){\comp}\iota_{V'}{\comp}\sigma(V')\end{array}}
	\arrow[2]{e,t,..}{\substack{P^{m}(\G{i}){\comp}\sigma(X){\comp}L\\ - P^{m}(\G{(i{\comp}L)}){\comp}\iota_{V'}{\comp}\sigma(V')}}
	\arrow[2]{s,=}
\node{}
\node{\begin{array}{c}P^{m}(\G{i}){\comp}\sigma(X){\comp}f{\comp}g\\ - P^{m}(\G{(i{\comp}f)}){\comp}\iota_V{\comp}\sigma(V){\comp}g\end{array}}
	\arrow[2]{s,=}
\\
\node{}
\node{(C)}
\node{}
\\
\node{\begin{array}{c}P^{m}(\G{(\tilde{L}{\comp}i')}){\comp}(\sigma(X){\comp}e^{X}_1 \\ - \iota_{X}){\comp}\Sigma\ad(f')\end{array}}
	\arrow[2]{e,t,..}{\substack{P^{m}(\G{i}){\comp}(\sigma(X){\comp}e^X_1 \\ - \iota_X){\comp}\Sigma\ads(L)}}
	\arrow[2]{s,..}
\node{}
\node{\begin{array}{c}P^{m}(\G{i}){\comp}(\sigma(X){\comp}e^X_1 \\ - \iota_X){\comp}\Sigma\ad(f{\comp}g)\end{array}}
	\arrow[2]{s,..}
\\
\node{}
\node{(D)}
\node{}
\\
\node{P^{m}(\G{(\tilde{L}{\comp}i')}){\comp}p^{\G{X}}_m{\comp}H^{\sigma(X)}_m(e^{X}_1){\comp}\Sigma\ad(f')}
	\arrow[2]{e,..}
	\arrow[2]{s,=}
\node{}
\node{P^{m}(\G{i}){\comp}p^{\G{X}}_m{\comp}H^{\sigma(X)}_m(e^X_1){\comp}\Sigma\ad(f{\comp}g)}
	\arrow[2]{s,=}
\\
\node{}
\node{(E)}
\node{}
\\
\node{P^{m}(\G{(\tilde{L}{\comp}i')}){\comp}p^{\G{X}}_m{\comp}H^{\sigma(X)}_m(f')}
	\arrow[2]{e,t}{\substack{P^{m}(\G{i}){\comp}p^{\G{X}}_m{\comp}\\H^{\sigma(X)}_m(e^X_1){\comp}\Sigma\ads(L)}}
\node{}
\node{P^{m}(\G{i}){\comp}p^{\G{X}}_m{\comp}H^{\sigma(X)}_m(f{\comp}g)}
\end{diagram}
\end{equation}
where $c(-)$ denotes a constant homotopy and $\iota_X : P^1(\G{X}) \hookrightarrow P^m(\G{X})$, $\iota_V : P^1(\G{V}) \hookrightarrow P^m(\G{V})$ and $\iota_{V'} : P^1(\G{V'}) \hookrightarrow P^m(\G{V'})$ are inclusions.
The homotopy commutativities of $(C)$ and $(E)$ are trivial.
\par
To show the homotopy commutativity of $(A)$, we define a map $\theta_A : [0,1]{\times}[0,1]{\times}V' \to W$.
\begin{equation*}
\theta_A(s,t,u{\wedge}v'_0) = 
\begin{cases}
   L(s,\left(\frac{2u}{2-t}\right){\wedge}v'_0) \in W,\quad u \leq \frac{2-t}{2}\\
   {\ast} \in X \subset W,\quad u \geq \frac{2-t}{2}
\end{cases}
\quad\text{for $s, t \in [0,1]$ and $u{\wedge}v'_0 \in V'$.}
\end{equation*}
We have $\theta_A(s,0,u{\wedge}v'_0)=L(s,u{\wedge}v'_0)$, $\theta_A(s,1,u{\wedge}v'_0)=(L-c(0))(s,u{\wedge}v'_0)$ and that $\theta_A(0,t,u{\wedge}v'_0)$ and $\theta_A(1,t,u{\wedge}v'_0)$ give canonical homotopies from $f'(u{\wedge}v'_0)$ and $f{\comp}g(u{\wedge}v'_0)=f(u{\wedge}g_{0}(v'_{0}))$ to $(f'-0)(u{\wedge}v'_0)$ and $(f{\comp}g-0)(u{\wedge}v'_0)=(f-0)(u{\wedge}g_{0}(v'_{0}))$.
Thus $\theta_A$ gives a homotopy
\begin{equation*}
\begin{diagram}
\node{f'}
	\arrow[2]{e,t,..}{L}
	\arrow{s,l,..}{}
\node{}
\node{f{\comp}g}
	\arrow{s,r,..}{}
\\
\node{f' - 0}
	\arrow[2]{e,b,..}{L - c(0)}
\node{}
\node{f{\comp}g - 0.}
\end{diagram}
\end{equation*}
where the upper and lower rows are given by $t=0$ and $t=1$ and the left and right columns are given by $s=0$ and $s=1$.
By applying the composition with $P^{m}(\G{i}){\comp}\sigma(X)$ from the left, we get the homotopy commutativity of $(A)$.
\par
To show the homotopy commutativity of $(B)$, we define a map $\theta_B : [0,1]{\times}[0,1]{\times}V' \to W$.
\begin{equation*}
\theta_B(s,t,v') = 
\begin{cases}
   \left(\frac{2t}{1+s}\right){\wedge}g(v') \in W,\quad t \leq \frac{1+s}{2}\\
   L(2+s-2t,v') \in X \subset W,\quad t \geq \frac{1+s}{2}
\end{cases}
\quad\text{for $s, t \in [0,1]$ and $v' \in V'$.}
\end{equation*}
We have $\theta_B(s,0,v')= \ast$, $\theta_B(0,t,v')=\tilde{L}{\comp}F'(t,v')$, $\theta_B(s,1,v')=i{\comp}L(s,v')$ and $\theta_B(1,t,v')=F(t,g(v'))$.
Thus $\theta_B$ gives a homotopy
\begin{equation*}
\begin{diagram}
\node{0}
	\arrow[5]{e,t,..}{c(0)}
	\arrow[2]{s,l,..}{P^{m}(\G(\tilde{L}{\comp}F'))}
\node{}
\node{}
\node{}
\node{}
\node{0}
	\arrow[2]{s,r,..}{P^{m}(\G{F})}
\\
\node{}
\\
\node{P^{m}(\G{(\tilde{L}{\comp}i'{\comp}f')})}
	\arrow[5]{e,b,..}{P^{m}(\G{(i{\comp}L)})}
\node{~~}
\node{~~}
\node{~~}
\node{~~}
\node{P^{m}(\G{(i{\comp}f{\comp}g)})}
\end{diagram}
\end{equation*}
where the upper and lower rows are given by $t=0$ and $t=1$ and the left and right columns are given by $s=0$ and $s=1$.
By applying the composition with $\iota_{V'}{\comp}\sigma(V')$ from the right, we get the homotopy commutativity of $(B)$.
\par
To show the homotopy commutativity of $(D)$, we fix a homotopy 
\begin{equation*}
\sigma(X){\comp}e^{X}_1 - \iota_X \sim p^{\G{X}}_m{\comp}H^{\sigma(X)}_m(e^{X}_1)
\end{equation*}
as a map $K : [0,1]{\times}{\Sigma}\G{X} \to P^m(\G{X})$ with $K(0,-)=(\sigma(X){\comp}e^{X}_1 - \iota_X)(-)$ and $K(1,-)=p^{\G{X}}_m{\comp}H^{\sigma(X)}_m(e^{X}_1)(-)$.
Using it, we define a map $\theta_D : [0,1]{\times}[0,1]{\times}V' \to P^m(\G{X})$.
\begin{equation*}
\theta_D(s,t,u{\wedge}v'_0) = K(t,u{\wedge}\ad(L)(s,v'_0))
\quad\text{for $s, t \in [0,1]$ and $u{\wedge}v'_0 \in V'$.}
\end{equation*}
We have $\theta_D(s,0,u{\wedge}v'_0)=(\sigma(X){\comp}e^{X}_1 - \iota_X)(u{\wedge}\ad(L)(s,v'_0))$, $\theta_D(0,t,u{\wedge}v'_0)=K(t,u{\wedge}\ad(f')(v'_0))$, $\theta_D(s,1,u{\wedge}v'_0)=p^{\G{X}}_m{\comp}H^{\sigma(X)}_m(e^{X}_1)(u{\wedge}\ad(L)(s,v'_0))$ and $\theta_D(1,t,u{\wedge}v'_0)=K(t,u{\wedge}\ad(f{\comp}g)(v'_0))$.
Thus $\theta_D$ gives a homotopy
\begin{equation*}
\dgHORIZPAD=2em
\begin{diagram}
\node{(\sigma(X){\comp}e^{X}_1 - \iota_X){\comp}{\Sigma}\ad(f')}
	\arrow[5]{e,t,..}{(\sigma(X){\comp}e^{X}_1 - \iota_X){\comp}{\Sigma}\ad(L)}
	\arrow{s,l,..}{}
\node{~}
\node{~}
\node{~}
\node{~}
\node{(\sigma(X){\comp}e^{X}_1 - \iota_X){\comp}{\Sigma}\ad(f{\comp}g)}
	\arrow{s,r,..}{}
\\
\node{p^{\G{X}}_m{\comp}H^{\sigma(X)}_m(e^{X}_1){\comp}{\Sigma}\ad(f')}
	\arrow[5]{e,b,..}{p^{\G{X}}_m{\comp}H^{\sigma(X)}_m(e^X_1){\comp}\Sigma\ads(L)}
\node{~}
\node{~}
\node{~}
\node{~}
\node{p^{\G{X}}_m{\comp}H^{\sigma(X)}_m(e^{X}_1){\comp}{\Sigma}\ad(f{\comp}g)}
\end{diagram}
\end{equation*}
where the upper and lower rows are given by $t=0$ and $t=1$ and the left and right columns are given by $s=0$ and $s=1$.
By applying the composition with $P^{m}(\G{i})$ from the left, we get the homotopy commutativity of $(D)$.
\par
The maps $\sigma'(W){\comp}\tilde{L}$ and $P^{m+1}(\G(\tilde{L})){\comp}\sigma'(W')$ are given by the homotopies in the diagram (\ref{eq:homotopy-of-homotopy}) indicated by the top and right arrows and left and bottom arrows respectively, and hence the homotopy commutativity of the homotopies in the diagram (\ref{eq:homotopy-of-homotopy}) implies that $\sigma'(W){\comp}\tilde{L} \sim P^{m+1}(\G{\tilde{L}}){\comp}\sigma'(W')$.
\end{Proof}
\begin{Rem}\label{rem:category-W}
Since $e^{W}_{m+1}{\comp}\sigma'(W)$ and the identity $1_{W}$ coincide on $X$ up to homotopy, the difference between them with respect to the co-action of $\Sigma{V}$ is given by a map $\gamma : \Sigma{V} \to W$.
Here, we know the fibration (\ref{eq:Stasheff-fibration}), induces the following split short exact sequence:
\begin{equation}\label{eq:Sugawara-exact}
0 \to [\Sigma{V},E^{m+1}(\G{W})] \overset{{p^{\G{W}}_m}_{\ast}}\to [\Sigma{V},P^{m}(\G{W})] \overset{{e^W_m}_{\ast}}\to [\Sigma{V},W] \to 0.
\end{equation}
Thus $\gamma$ can be pulled back to a map $\gamma_0 : \Sigma{V} \to P^{m}(\G{W}) \subset P^{m+1}(\G{W})$.
By adding $\gamma_0$ to $\sigma'(W)$, we obtain a genuine compression $\sigma(W)$ of $1_{W}$.
We often call this $\sigma(W)$ a `standard' structure map for $\cat{W} \leq m+1$ which gives a subset of $H^S_m(f)$ for $f \in [V',W]$ as follows:
\begin{equation*}
H^{SS}_m(f) = \{H^{\sigma(W)}_m(f) \,\vert\, \text{$\sigma(W)$ is a `standard' structure map for $\cat{W} \leq m+1$}\}
\end{equation*}
\end{Rem}

The cofibration sequence $V \overset{f}\to X \overset{i}\to W \overset{q}\to {\Sigma}V$ induces the following cofibration sequence:
\begin{equation*}
V{\ast}S^{n-1} \overset{\hat{f}}\to W \cup X{\times}S^n \overset{\hat{i}}\to W{\times}S^n \overset{\hat{q}}\to {\Sigma}^{n+1}V;
\end{equation*}
here $\hat{f}$ is given by the relative Whitehead product $[\chi_{f},1_{S^{n}}]^r : V{\ast}S^{n-1} \to W \cup X{\times}S^n$, where $\chi_f : (C(V),V) \to (W,X)$ denotes the characteristic map given as a relative homeomorphism.

By the proof of Proposition 5.8 in \cite{Iwase:counter-ls}, the following result is obtained using Remark \ref{rem:category-W'}.

\begin{Prop}\label{prop:obstruction2}
For a given structure map $\sigma(X)$ for $\cat{X} = m$, the map $\hat{f}$ makes the following diagram without the dotted arrows commute up to homotopy.
\begin{equation*}
\begin{diagram}
\divide\dgARROWLENGTH by2
\dgHORIZPAD=.5em
\dgVERTPAD=1ex
\node{V{\ast}S^{n-1}}
        \arrow{e,t}{\hat{f}}
        \arrow{s,r}{H_{m}(f){\ast}1_{S^{n-1}}}
\node{W \cup X{\times}S^n}
        \arrow{e,t,J}{\hat{i}}
        \arrow{s,r}{\sigma'(W){\times}1_{S^n}\vert_{W \cup X{\times}S^n}}
\node{W{\times}S^n}
        \arrow[2]{s,r,..}{\lower-6.5ex\hbox{${}_{\sigma(W{\times}S^n)}$}}
        \arrow{e,t}{\hat{q}}
\node{{\Sigma}^{n+1}V}
\\
\node{E^{m+1}(\G{X}){\ast}{S^{n-1}}}
        \arrow{s,r,L}{E^{m+1}(\G{i}){\ast}j_{n-1}}
        \arrow{e,b,..}{\hat{p}^{\G{X}}_{m}}
\node{P^{m+1}(\G{W}) \cup P^{m}(\G{X}){\times}S^n}
        \arrow{s,r,L}{}
\\
\node{E^{m+1}(\G{W}){\ast}\G{S^n}}
        \arrow{e,b}{\hat{p}^{\G{W}}_{m}}
\node{P^{m+1}(\G{W}) \cup P^{m}(\G{W}){\times}S^n}
        \arrow{e,b,J}{}
\node{P^{m+1}(\G{W}){\times}S^n}
        \arrow{e,b}{e^{W}_{m+1}{\times}1_{S^n}}
\node{W{\times}S^n,}
\end{diagram}
\end{equation*}
where $\hat{p}^{\G{X}}_{m} = [\chi_{p^{\G{X}}_{m}},1_{S^{n}}]^r$, $\hat{p}^{\G{W}}_{m} = [\chi_{p^{\G{W}}_{m}},e^{S^{n}}_1]^r$ and $j_{n-1}$ denotes the bottom cell inclusion $S^{n-1} \subset \G{S^n}$ and $\sigma'(W)$ is the extension of $\sigma(X)$ by Proposition \ref{prop:obstruction} and Remark \ref{rem:category-W'}.
\end{Prop}

We then have the following result.

\begin{Thm}\label{thm:stably-trivial}
The following two statements hold for $W$ with $\cat{W} \leq m+1$, $m = \cat{X}$.
\begin{description}
\item{(1)}\label{thm:stably-trivial-1}
$\cat{W} \leq \cat{X}$ if the set $H^{S}_m(f)$ contains the trivial element.
\item{(2)}\label{thm:stably-trivial-2}
$\cat{W{\times}S^{n}} \leq \cat{X}+1$ if the set ${\Sigma}^{n}_{\ast}H^{S}_m(f)$ contains the trivial element.
\end{description}
\end{Thm}
\begin{Cor}\label{cor:stably-trivial}
Let $W$ be the space constructed as in the above theorem with $\cat{W}=m+1$.
Then $W$ is a counter example to Ganea's conjecture if the set ${\mathcal H}^S_m(f)$ contains the trivial element.
\end{Cor}
\begin{Proof*}{\it Proof of Theorem \ref{thm:stably-trivial}.}
Firstly we show (1):
If $E^{m+1}(\G{i}){\comp}H^{\sigma(X)}_m(f) : V \to E^{m+1}(\G{W})$ is trivial for some structure map $\sigma(X)$, then, by Proposition \ref{prop:obstruction} and Remark \ref{rem:category-W'}, $\sigma(X)$ is extendible to $W$.
By the argument given in Remark \ref{rem:category-W}, we obtain a genuine compression $\sigma(W) : W \to P^m(\G{W})$ of $1_W$.
Thus $\cat{W} \leq m$.

Secondly we show (2):
If $(E^{m+1}(\G{i}){\ast}j_{n-1}){\comp}{\Sigma}^{n}H^{\sigma(X)}_m(f)$ $:$ $V{\ast}S^{n-1}$ $\to$ $E^{m+1}(\G{W}){\ast}\G{S^{n}}$ is trivial for some structure map $\sigma(X)$, then, by Proposition \ref{prop:obstruction2}, the map $\sigma'(W){\times}1_{S^n}\vert_{W \cup X{\times}S^n} : W \cup X{\times}S^n \to P^{m+1}(\G{W}) \cup P^{m}(\G{X}){\times}S^n \subset P^{m+1}(\G{W}) \cup P^{m}(\G{W}){\times}S^n$ is extendible to $W{\times}S^{n}$.
By the argument given in Remark \ref{rem:category-W} together with the fact that the natural map $[\Sigma^{n+1}{V},Z{\times}B{\cup}C{\times}Y]$ $\to$ $[\Sigma^{n+1}{V},Z{\times}Y]$ $\cong$ $[\Sigma^{n+1}{V},Z]{\times}[\Sigma^{n+1}{V},Y]$ is split surjective for any pointed pairs $(Z,C)$ and $(Y,B)$, we obtain a compression $\sigma'' : W{\times}S^{n} \to P^{m+1}(\G{W}) \cup P^m(\G{W}){\times}S^{n}$ of $\sigma'(W){\times}1_{S^n}$, relative to $\sigma'(W){\times}1_{S^n}\vert_{W \cup X{\times}S^n}$.
Also by Remark \ref{rem:category-W}, the difference between $e^{W}_{m+1}{\comp}\sigma'(W)$ and $1_{W}$ with respect to the co-action of $\Sigma{V}$ can be pulled back to $\gamma_0 :\Sigma{V} \to {\Sigma}\G{W} \subset P^m(\G{W})$.
%
%
We define a map $\sigma(W{\times}S^n)$ by 
\begin{equation*}
\multiply\dgARROWLENGTH by2
\begin{diagram}
\sigma(W{\times}S^n) : W{\times}S^n \overset{\mu{\times}1_{S^n}}\to %
\node{(W \vee \Sigma{V}){\times}S^n}
\arrow{e,t}{\sigma'' \cup (\gamma_0{\times}1_{S^n})} \node{P^{m+1}({\Omega}W){\times}\{\ast\} \cup P^m(\Omega{W}){\times}S^n,}
\end{diagram}
\end{equation*}
where $\mu : W \to W \vee {\Sigma}V$ denotes the co-action.
Since $\sigma''$ is homotopic to $\sigma'(W){\times}1_{S^n}$ in $P^{m+1}(\Omega{W}){\times}S^n$ relative to the subspace $\{\ast\}{\times}S^n$, $\sigma(W{\times}S^n)$ is homotopic to $(\sigma'(W)+\gamma){\times}1_{S^n}$ in $P^{m+1}(\Omega{W}){\times}S^n$, which is a compression of $1_{W}{\times}1_{S^n}$.
Thus $\cat{W{\times}S^{n}} \leq m+1$, by Fact \ref{cor:product-formula-sphere}.
\end{Proof*}

Let $V$ be a ($e-1$)-connected co-H-space and $X$ a ($d-1$)-connected CW complex, $e \geq d \geq 2$.

\begin{Thm}\label{thm:stably-non-trivial}
Let $n \geq 1$ and suppose $X$ to satisfy $\dim{X} \leq d{\cdot}m + d - 2$, $m = \cat{X}$.
Then the following two statements hold for $W$ and $H_m(f)=H^{\sigma(X)}_m(f)$ where $\sigma(X)$ is the unique structure map, by Proposition \ref{prop:uniqueness}.
\begin{description}
\item{(1)}
$\cat{W} = m + 1$ if $E^{m+1}(\G{i}){\comp}H_m(f) \ne 0$.
\item{(2)}
$\cat{W{\times}S^{n}} = m + 2$ if $(E^{m+1}(\G{i}){\ast}j_{n-1}){\comp}{\Sigma}^{n}H_m(f) \ne 0$ for sufficiently large $n$.
\item{(3)}
$\cat{W{\times}S^{n}} = m + 2$ if $(E^{m+1}(\G{i}){\ast}j_{n}){\comp}{\Sigma}^{n+1}H_m(f) \ne 0$.
\end{description}
\end{Thm}
\begin{Rem}\label{rem:remove-inc}
If the condition $\dim{V} < d\cdot\cat{X} + e - 1$ is satisfied, then for dimensional reasons, we may remove the maps `$E^{m+1}(\G{i})$', `$E^{m+1}(\G{i}){\ast}j_{n-1}$' and `$E^{m+1}(\G{i}){\ast}j_{n}$' from the statements in Theorem \ref{thm:stably-non-trivial}.
\end{Rem}
\begin{Cor}\label{cor:counter-ls}
Let $W$ be the space constructed as the mapping cone of $f : V \to X$ from simply connected spaces $V$ and $X$ with $\dim{X} \leq d\cdot\cat{X} + d - 2$ and $\dim{V} \leq d\cdot\cat{X} + e - 2$, where $V$ and $X$ are ($e-1$)-connected and ($d-1$)-connected, resp.
If $\cat{W{\times}S^k} = \cat{W}$ for some $k$, then $\cat{W{\times}S^n} = \cat{W}$ for all $n \geq k$.
Thus Conjecture \ref{conj:iwase-number} is true for such $W$.
\end{Cor}
\begin{Cor}\label{cor:criterion-ls}
Let $W$ be the space constructed as the mapping cone of $f : V \to X$ from simply connected spaces $V$ and $X$ with $\dim{X} \leq d\cdot\cat{X} + d - 2$ and $\dim{V} \leq d\cdot\cat{X} + e - 2$, where $V$ and $X$ are ($e-1$)-connected and ($d-1$)-connected, resp.
Then Ganea's conjecture for $W$ is true if and only if the unique stabilised higher Hopf invariant ${\mathcal H}_m(f)$ is non-trivial.
\end{Cor}
\begin{Proof*}{\it Proof of Theorem \ref{thm:stably-non-trivial}.}
Let $\dim{X} = k$.
By the assumption, $\cat{W} \leq m+1$.

Firstly we show (1):
If $\cat{W} \leq m$, then there exists a compression $\sigma$ $:$ $W$ $\to$ $P^{m}(\G{W})$ of the identity $1_{W}$.
Since the pair $(P^{\infty}(\G{W}),P^{m}(\G{W}))$ is ($d(m+1)-1$)-connected, the inclusion map $\iota_{m} : P^{m}(\G{W}) \to P^{\infty}(\G{W})$ induces a bijection of homotopy sets $[Z,P^{m}(\G{W})] \to [Z,P^{\infty}(\G{W})]$ for any space $Z$ of dimension up to $d(m+1)-2$.
Then, for dimensional reasons, it follows that the restriction $\sigma_0 = \sigma\vert_{X}$ is unique up to homotopy in $P^{m}(\G{W})$.
Hence we may assume that $\sigma_0$ equals to $\sigma(X)$, the unique structure map for $\cat{X} = m$ by Proposition \ref{prop:uniqueness}.
Hence $P^m(\G{i}){\comp}\sigma(X){\comp}f \sim p^{\G{W}}_{m}{\comp}E^{m+1}(\G{i}){\comp}H_m(f)$ is null-homotopic.
As $P^m(\G{W})_{\ast}$ in (\ref{eq:Sugawara-exact}) is a monomorphism for $m \geq 1$, $E^{m+1}(\G{i}){\comp}H_m(f)$ is null-homotopic.
Thus the existence of the compression $\sigma$ of $1_W$ implies the triviality of $E^{m+1}(\G{i}){\comp}H_m(f) : V \to E^{m+1}(\G{W})$.

Secondly we show (2):
Let $\sigma(W)$ be a structure map for $\cat{W}{\leq}m+1$ obtained from $\sigma(X)$ by Proposition \ref{prop:obstruction}, Remarks \ref{rem:category-W'} and \ref{rem:category-W}, and let $\hat{j}_{m,n} : P^{m+1}(\G{W}) \cup P^{m}(\G{W}){\times}S^{n} \subset P^{\infty}(\G{W}){\cup}P^{m}(\G{W}){\times}S^{n}$ and $\hat{i}_{m,n} : P^{\infty}(\G{W}) \cup P^{m}(\G{W}){\times}S^{n} \subset P^{\infty}(\G{W}){\times}S^{n}$ be inclusions.
If $\cat{W{\times}S^{n}} \leq m+1$, then by Fact \ref{cor:product-formula-sphere}, there exists a compression $\sigma : W{\times}S^{n} \to P^{m+1}(\G{W}) \cup P^{m}(\G{W}){\times}S^{n}$ of the identity $1_{W{\times}S^{n}}$.
Hence we have 
\begin{align*}&
(e^{W}_{\infty}{\times}1_{S^n}){\vert}_{P^{\infty}(\G{W}){\cup}P^{m}(\G{W}){\times}S^{n}}{\comp}\hat{j}_{m,n}{\comp}\sigma{\vert}_{W{\cup}X{\times}S^{n}} \sim 1_{W{\times}S^{n}}{\vert}_{W{\cup}X{\times}S^{n}} \quad\text{and}\quad
\\&
(e^{W}_{\infty}{\times}1_{S^n}){\vert}_{P^{\infty}(\G{W}){\cup}P^{m}(\G{W}){\times}S^{n}}{\comp}\hat{j}_{m,n}{\comp}(\sigma(W){\times}1_{S^{n}}{\vert}_{W{\cup}X{\times}S^{n}}) \sim 1_{W{\times}S^{n}}{\vert}_{W{\cup}X{\times}S^{n}}.
\end{align*}
Here we know the pair $(P^{\infty}(\G{W}){\times}S^{n},P^{\infty}(\G{W}){\cup}P^{m}(\G{W}){\times}S^{n})$ is ($dm+n+d-1$)-connected and there is a co-action $\mu'_{W} : W{\cup}X{\times}S^n \to W{\cup}X{\times}S^n \vee {\Sigma}V$ associated with the cofibration sequence $V \to X{\times}S^n \to W \cup X{\times}S^n$.
Since the subspace $X{\times}S^{n}$ is of dimension $k+n \leq dm+d+n-2$, the restrictions $\hat{j}_{m,n}{\comp}\sigma\vert_{X{\times}S^n}$ and $\hat{j}_{m,n}{\comp}\sigma(W){\times}1_{S^n}\vert_{X{\times}S^n} = \sigma(X){\times}1_{S^n}$ are homotopic and the difference between $\hat{j}_{m,n}{\comp}\sigma$ and $\sigma(W){\times}1_{S^n}\vert_{W{\cup}X{\times}S^n}$ with respect to the co-action $\mu'_{W}$ is given by a map into $E^{m+1}(\G{W}){\ast}\G{S^{n}}$ the fibre of $e^{W}_{\infty}{\times}1_{S^n}{\vert}_{P^{\infty}(\G{W}){\cup}P^{m}(\G{W}){\times}S^{n}}$:
\begin{equation*}
\sigma{\vert}_{W \cup X{\times}S^{n}} \sim \sigma(W){\vert}_{W \cup X{\times}S^{n}} + \gamma, \qquad \gamma=\hat{j}_{m,n}{\comp}\hat{p}^{\G{W}}_m{\comp}\gamma_0,\quad \gamma_0 : {\Sigma}V \to E^{m+1}(\G{W}){\ast}\G{S^{n}}.
\end{equation*}
When $n > \dim{V}-d(m+1)+2$, we also obtain that $\hat{j}_{m,n}{\comp}\sigma\vert_{W{\cup}X{\times}S^n} \sim \hat{j}_{m,n}{\comp}\sigma(W){\times}1_{S^n}\vert_{W{\cup}X{\times}S^n}$ for dimensional reasons.
In this case, we have 
\begin{align*}&
\hat{j}_{m,n}{\comp}\hat{p}^{\G{W}}_m{\comp}(E^{m+1}(\G{i}){\ast}j_{n-1}){\comp}{\Sigma}^{n}H_m(f)
\sim \hat{j}_{m,n}{\comp}(\sigma(W){\times}1_{S^{n}}{\vert}_{W{\cup}X{\times}S^{n}}){\comp}\hat{f}
\sim \hat{j}_{m,n}{\comp}\sigma{\vert}_{W{\cup}X{\times}S^{n}}{\comp}\hat{f}
\end{align*}
which is trivial, since $\sigma{\vert}_{W{\cup}X{\times}S^{n}}$ is extendible on $W{\times}S^{n}$.
Thus $\cat{W{\times}S^{n}} \leq m+1$ implies $(E^{m+1}(\G{i}){\ast}j_{n-1}){\comp}{\Sigma}^{n}H_m(f) \sim 0$, provided that $n > \dim{V}-d(m+1)+2$.

Thirdly we show (3):
In this case, the difference $\gamma$ might not be trivial evenif $\cat{W{\times}S^{n}} \leq m+1$.
Since the composition $\proj_1{\comp}\hat{j}_{m,n}{\comp}\sigma : W{\times}S^{n} \overset{\hat{j}_{m,n}{\comp}\sigma}\to P^{\infty}(\G{W}) \cup P^{m}(\G{W}){\times}S^{n} \overset{\proj_1}\to P^{\infty}(\G{W})$ factors through $P^{\infty}(\G{W}){\times}S^n$ as $\proj_1{\comp}\hat{j}_{m,n}{\comp}\sigma = \proj_1{\comp}\hat{i}_{m,n}{\comp}\hat{j}_{m,n}{\comp}\sigma$, we have $e^{W}_{\infty}{\comp}\proj_1{\comp}\hat{j}_{m,n}{\comp}\sigma = \proj_1{\comp}(e^{W}_{\infty}{\times}1_{S^n}){\comp}\hat{i}_{m,n}{\comp}\hat{j}_{m,n}{\comp}\sigma \sim \proj_1$, which makes the following diagram commutative up to homotopy.
\begin{equation*}
\divide\dgARROWLENGTH by2
\dgHORIZPAD=.5em
\dgVERTPAD=1ex
\begin{diagram}
\node{}
\node{W}
	\arrow[2]{e,t}{(e^{W}_{\infty})^{-1}}
\node{}
\node{P^{\infty}(\G{W})}
\\
\node{W {\cup} X{\times}S^{n}}
	\arrow{e,J}
\node{W{\times}S^{n}}
	\arrow{n,l}{\proj_1}
	\arrow[2]{e,t}{\hat{j}_{m,n}{\comp}\sigma}
	\arrow{s,l}{\proj_1}
\node{}
	\arrow{e,!}
	\arrow{w,!}
	\arrow{s,!}
\node{P^{\infty}(\G{W}){\cup}P^{m}(\G{W}){\times}S^{n}}
	\arrow{n,r}{\proj_1}
	\arrow{s,r}{\proj_1}
\\
\node{}
\node{W}
	\arrow[2]{e,t}{(e^{W}_{\infty})^{-1}}
\node{}
\node{P^{\infty}(\G{W}).}
\end{diagram}
\end{equation*}
Taking push-outs of both right and center columns of the above diagram, we have a map 
\begin{equation*}
\sigma_1 : W{\times}S^{n+1} \to P^{\infty}(\G{W}){\cup}P^{m}(\G{W}){\times}S^{n+1}\quad\text{with}\quad \sigma_1{\comp}(1_{W}{\times}E_{n}) = {E}_{m,n}{\comp}\hat{j}_{m,n}{\comp}\sigma,
\end{equation*}
where $E_n$ is the inclusion $S^{n}=S^{n}{\times}\{0\} \hookrightarrow S^{n+1} = S^{n}{\times}[-1,1]/\text{$\sim$}$, $(x,-1) \sim (x,1) \sim (\ast,t)$ for $x \in S^n$ and $t \in [-1,1]$ and ${E}_{m,n}=1_{P^{\infty}(\G{W})}{\times}E_n{\vert}_{P^{\infty}(\G{W}){\cup}P^{m}(\G{W}){\times}S^{n}}$.
By the definition of $\sigma_1$, we have $\sigma_1{\vert}_{X{\times}S^{n+1}} \sim \hat{j}_{m,n+1}{\comp}(\sigma(X){\times}1_{S^{n+1}})$ and the difference between $\sigma_1{\vert}_{W{\cup}X{\times}S^{n+1}}$ and $\hat{j}_{m,n+1}{\comp}(\sigma(W){\times}1_{S^{n+1}}){\vert}_{W{\cup}X{\times}S^{n+1}}$ is given by ${E}_{m,n}{\comp}\gamma$.
Since $P^{\infty}(\G{W}){\cup}P^{m}(\G{W}){\times}(S^{n}{\times}0)$ is compressible into $P^{\infty}(\G{W})$, so is ${E}_{m,n}{\comp}{\gamma}$.
Thus ${E}_{m,n}{\comp}{\gamma} \sim \incl_1{\comp}\proj_1{\comp}{E}_{m,n}{\comp}{\gamma} = \incl_1{\comp}\proj_1{\comp}\gamma \sim \incl_1{\comp}\proj_1{\comp}\hat{i}_{m,n}{\comp}\hat{j}_{m,n}{\comp}\gamma \sim 0$, where $\incl_1$ or $\proj_1$ denotes an appropriate inclusion or projection.
Then by using the same argument as in (2), we obtain $(E^{m+1}(\G{i}){\ast}j_{n}){\comp}{\Sigma}^{n+1}H_m(f) \sim 0$, if $\cat{W{\times}S^{n}} \leq m+1$.
\end{Proof*}

\section{Higher Hopf invariants and the reduced diagonal}\label{sect:hopf-diag}

We state here another property of the higher Hopf invariants, which is a generalisation of Theorem 5.14 of Boardman-Steer \cite{BS:hopf-invariants}.

\begin{Thm}\label{thm:reduced-diagonal}
Let $V \overset{f}\to X \overset{i}\to W$ be a cofibration sequence as above with $\cat{X} = m$.
Then the reduced diagonal $\hat{\Delta}_{m+1} : W \to W^{m+1} \to \wedge^{m+1}W$ denotes homotopic to the composition $\wedge^{m+1}e^W_{1}{\comp}{\Sigma}h^W_{m+1}{\comp}{\Sigma}(E^{m+1}(\G{i}){\comp}H_{m}(f)){\comp}q$, where $q$ is the collapsing map $W \to W/X={\Sigma}V$.
\end{Thm}

\begin{Cor}\label{cor:reduced-diagonal}
Let $V$, $X$, $f$ and $W$ be as above with $\cat{X} = m$ and $\cat{W} = m+1$.
Then the $n$-fold suspension of the reduced diagonal $\hat{\Delta}_{m+1}{\wedge}1_{S^n} : W{\times}S^n \to W^{m+1}{\times}S^n \to \wedge^{m+1}W{\wedge}S^n$ is homotopic to the composition ${\Sigma}^{n}(\wedge^{m+1}e^W_{1}{\comp}{\Sigma}h^W_{m+1}{\comp}{\Sigma}(E^{m+1}(\G{i}){\comp}H_{m}(f))){\comp}\hat{q}$, where $\hat{q}$ denotes the collapsing map $W{\times}S^n \to W{\times}S^n/(W \cup X{\times}S^n) = (W/X){\wedge}S^n = {\Sigma}^{n+1}V$.
\end{Cor}

For any space $Z$, Ganea showed that there is a commutative ladder of fibrations up to homotopy (see [5] or the proof of Theorem 1.1 of \cite{Iwase:counter-ls}):
\begin{equation}\label{eq:pb-diagram}
\dgHORIZPAD=.5em
\dgVERTPAD=1ex
\begin{diagram}
\node{{E}^{m+1}(\G{Z})}
	\arrow{e,t}{p^Z_{m}}
	\arrow{s,r}{h^Z_{m+1}}
\node{P^{m}(\G{Z})}
	\arrow{e,t}{e^Z_{m}}
	\arrow{s,r}{\Delta'_{m+1}}
\node{Z}
	\arrow{s,r}{\Delta_{m+1}}
\\
\node{\bar{E}^{m+1}(\G{Z})}
	\arrow{e,t}{q^Z_{m+1}}
\node{Z^{[m+1]}}
	\arrow{e,t,J}{}
\node{Z^{m+1},}
\end{diagram}
\end{equation}
where 
$q^Z_{m+1} : \bar{E}^{m+1}(\G{Z}) \to Z^{[m+1]}$ is given by $q^Z_{m+1}(\Sigma t_i\ell_i)$ $=$ $(\ell_0(t_0/t_M),...,\ell_m(t_m/t_M))$, $t_M$ $=$ $\max(t_0,...,t_{m})$ and 
$\Delta'_{m+1}$ denotes a map which makes the right hand square of (\ref{eq:pb-diagram}) a homotopy pull-back diagram.

Using the precise description of the fibration $q^Z_{m+1}$, we obtain the following proposition.
\begin{Prop}\label{prop:diagonal}
There exists a map $e_{m+1} : Z^{[m+1]} \cup \hat{C}(\bar{E}^{m+1}(\G{Z})) \to Z^{m+1}$ which makes the following diagram commute up to homotopy, where $\hat{C}$ denotes the functor taking unreduced cones with base point $0{\wedge}*$:
\begin{equation*}
\dgHORIZPAD=.5em
\begin{diagram}
\node{\bar{E}^{m+1}(\G{Z})}
        \arrow{e,t}{q^{Z}_{m+1}}
        \arrow{s,=}
\node{Z^{[m+1]}}
        \arrow{e,t,J}{}
        \arrow{s,=}
\node{Z^{[m+1]} \cup \hat{C}(\bar{E}^{m+1}(\G{Z}))}
        \arrow{e,t}{}
        \arrow{s,r}{e_{m+1}}
\node{\wedge^{m+1}{\Sigma}\G{Z}}
        \arrow{s,r}{\wedge^{m+1}e^{Z}_{1}}
\\
\node{\bar{E}^{m+1}(\G{Z})}
        \arrow{e,t}{q^{Z}_{m+1}}
\node{Z^{[m+1]}}
        \arrow{e,t,J}{}
\node{Z^{m+1}}
        \arrow{e,t}{}
\node{\wedge^{m+1}Z,}
\end{diagram}
\end{equation*}
where 
the upper row is a cofibration sequence.
\end{Prop}
\begin{Proof}
We define $e_{m+1} : Z^{[m+1]} \cup \hat{C}(\bar{E}^{m+1}(\G{Z})) \to Z^{m+1}$ by 
\begin{equation*}
e_{m+1}(t,\Sigma t_i\ell_i) = (\ell_0(t{\cdot}t_0/t_M),\cdots,\ell_m(t{\cdot}t_m/t_M)),
\end{equation*}
where $\sum_{0{\leq}i{\leq}m}t_i=1$ and $t_M = \max(t_0,...,t_m) \geq \frac{1}{m+1}$.
By collapsing $Z^{[m+1]}$, we obtain a map $\hat{e}_{m+1} : \hat{\Sigma}\bar{E}^{m+1}(\G{Z}) \to \wedge^{m+1}Z$ from $e_{m+1}$, where $\hat{\Sigma}$ denotes the functor taking unreduced suspensions with base point $0{\wedge}*$.
$\hat{e}_{m+1}$ is given by $\hat{e}_{m+1}(t \wedge \Sigma t_i\ell_i) = \ell_0(t t_0/t_M){\wedge}...{\wedge}\ell_m(t t_m/t_M)$.

Let us recall the homotopy equivalence $s_{m+1}$ $:$ $\hat{\Sigma}\bar{E}^{m+1}(\G{Z})$ $\to$ $\wedge^{m+1}\hat{\Sigma}\G{Z}$ given by
\begin{equation*}
s_{m+1}(t,\Sigma t_i\ell_i) = (u_0\wedge\ell_0,...,u_m\wedge\ell_m),\quad u_i = t_i\frac{t}{t_M},\ 0 \leq i \leq m.
\end{equation*}
Then we have the following relation of maps: $\wedge^{m+1}e^Z_1{\comp}h{\comp}s_{m+1} = \hat{e}_{m+1}$, where $h : \wedge^{m+1}\hat{\Sigma}\G{Z} \to \wedge^{m+1}{\Sigma}\G{Z}$ is the canonical homotopy equivalence.
Thus the proposition follows.
\end{Proof}

\begin{Proof*}{\it Proof of Theorem \ref{thm:reduced-diagonal}.}
By Remarks \ref{rem:category-W'} and \ref{rem:category-W} and Propositions \ref{prop:obstruction} and \ref{prop:diagonal}, we have the following commutative diagram up to homotopy.

\begin{equation*}
\begin{diagram}
\dgHORIZPAD=.5em
\dgVERTPAD=1ex
\node{V}
        \arrow{e,t}{f}
        \arrow{s,l}{E^{m+1}(\G{i}){\comp}H_{m}(f)}
\node{X}
        \arrow{e,t,J}{i}
        \arrow{s,r}{\sigma(X)}
\node{W}
        \arrow{e,t}{q}
        \arrow{s,r}{\sigma(W)}
\node{{\Sigma}V}
        \arrow{s,r}{{\Sigma}(E^{m+1}(\G{i}){\comp}H_{m}(f))}
\\
\node{E^{m+1}(\G{W})}
        \arrow{e,b}{}
        \arrow{s,l}{h^W_{m+1}}
\node{P^{m}(\G{W})}
        \arrow{e,b,J}{}
        \arrow{s,r}{}
\node{P^{m+1}(\G{W})}
        \arrow{e,b}{}
        \arrow{s,r}{\Delta_{m+1}{\comp}e^{W}_{m+1}}
\node{\Sigma{E^{m+1}(\G{W})}}
        \arrow{s,r}{\wedge^{m+1}e^W_{1}{\comp}{\Sigma}h^W_{m+1}}
\\
\node{\bar{E}^{m+1}(\G{W})}
        \arrow{e,b}{q^{W}_{m}}
\node{W^{[m+1]}}
        \arrow{e,b,J}{}
\node{W^{m+1}}
        \arrow{e,b}{\wedge}
\node{\wedge^{m+1}W.}
\end{diagram}%
\end{equation*}%
\end{Proof*}

\begin{Rem}\label{rem:reduced-diagonal}\quad
For any structure map $\sigma(X)$, by naturality, we have 
\begin{align*}&
\hat{\Delta}_{m+1} = (\wedge^{m+1}e^W_{1}){\comp}{\Sigma}h^W_{m+1}{\comp}{\Sigma}(E^{m+1}(\G{i}){\comp}H^{\sigma(X)}_{m}(f)){\comp}q
\\&\qquad 
 = (\wedge^{m+1}e^W_{1}){\comp}{\Sigma}h^W_{m+1}{\comp}{\Sigma}E^{m+1}(\G{i}){\comp}{\Sigma}H^{\sigma(X)}_{m}(f){\comp}q
\\&\qquad 
= (\wedge^{m+1}e^W_{1}){\comp}(\wedge^{m+1}{\Sigma}\G{i}){\comp}{\Sigma}h^X_{m+1}{\comp}{\Sigma}(H^{\sigma(X)}_{m}(f)){\comp}q
\\&\qquad 
= (\wedge^{m+1}{i}){\comp}(\wedge^{m+1}e^X_{1}){\comp}{\Sigma}h^X_{m+1}{\comp}{\Sigma}(H^{\sigma(X)}_{m}(f)){\comp}q
= (\wedge^{m+1}{i}){\comp}\bar{H}^{\sigma(X)}_m(f){\comp}q,
\end{align*}
where the composition $\bar{H}^{\sigma(X)}_m(f)$ $=$ $(\wedge^{m+1}e^X_{1}){\comp}{\Sigma}h^X_{m+1}{\comp}{\Sigma}(H^{\sigma(X)}_{m}(f))$ is the generalised version of the Berstein-Hilton crude Hopf invariant (see Definition \ref{def:crude-hopf-invariant}).
Thus we have
\begin{equation*}
(\wedge^{m+1}{i})_{\ast}{\comp}q^{\ast}\bar{H}^{S}_m(f) = \{\hat{\Delta}_{m+1}\}.
\end{equation*}
\end{Rem}

\section{Homology decomposition and product spaces}\label{sect:hom-decomp}

In this section, we always assume that $X$ is a connected finite complex with a homology decomposition $\{X_t, f_t : S_t(X) \to X_{t}\}_{t \geq 1}$ of $X$, where $S_t(X)$ is the Moore space of type $(H_{t+1}(X),t)$ for $t \geq 1$.
By modifying the arguments given in Curjel \cite{Curjel:k'inv_co-hopf}, we obtain the following result.
\begin{Thm}\label{thm:hom-decomp}
The homology decomposition $\{X_t, f_t : S_t(X) \to X_{t}\}_{t \geq 1}$ of a simply connected space $X$ satisfies $0 = \cat{X_1} \leq \cat{X_2} \leq ... \leq \cat{X_{t-1}} \leq \cat{X_t} \leq ... \leq \cat{X}$.
If, in addition, $\cat{X_{k-1}} = \cat{X_{k}} = m$, then we can choose structure maps $\sigma_{k-1}(X) : X_{k-1} \to P^{m}(\G{X_{k-1}})$ and $\sigma_{k}(X) : X_{k} \to P^{m}(\G{X_{k}})$ for $\cat{X_{k-1}} = m$ and $\cat{X_{k}} = m$ to be compatible with each other, i.e. $\sigma_{k}(X)\vert_{X_{k}} \sim \sigma_{k-1}(X)$ in $P^{m}(\G{X_{k}})$.
\end{Thm}
\begin{Proof}
To prove the former part of the theorem, it is sufficient to show $\cat{X_{k}} \leq \cat{X_{k+1}}$ for $k \geq 1$.
So we may assume that $X = X_{k+1}$ and $\cat{X} = m$.
If $m=0$, then Theorem \ref{thm:hom-decomp} is clearly true, and hence we may assume that $m \geq 1$.
Then there is a homotopy section $\sigma(X) : X \to P^{m}(\G{X})$.
By induction on $k$, we show the existence of a compression $\sigma_k(X) : X_{k} \to P^{m}(\G{X_{k}})$ of $\sigma(X)\vert_{X_k} : X_k \to P^{m}(\G{X})$.
In the case $k = 1$, we have $X_1 = \{\ast\}$, and hence the existence of $\sigma_1(X)$ is clear.
In the general case $k > 1$, by the induction hypothesis, we have a compression $\sigma_{k-1}(X) : X_{k-1} \to P^{m}(\G{X_{k-1}})$ of $\sigma(X){\comp}i_{k-1} : X_{k-1} \to P^{m}(\G{X})$, where $i_{t} : X_{t} \hookrightarrow X$ denotes the canonical inclusion, $t \geq 1$; in addition, $e^{X_{k-1}}_m{\comp}\sigma_{k-1}(X) \sim 1_{X_{k-1}}$.
%
Let us consider the following commutative diagram:
\begin{equation}\label{eq:cat-structure}
\begin{diagram}
\divide\dgARROWLENGTH by3
\dgHORIZPAD=.5em
\dgVERTPAD=1ex
\node{}
\node{}
\node{}
\node{}
\node{E^{m+1}(\G{X_k})}
	\arrow[3]{s,l}{p^{\G{X_{k}}}_{m}}
	\arrow[3]{e,t,J}{E^{m+1}(\G{i_{k}})}
\node{}
\node{}
\node{E^{m+1}(\G{X})}
	\arrow[3]{s,r}{p^{\G{X}}_{m}}
\\
\\
\node{X_{k-1}}
	\arrow{se,b}{\sigma_{k-1}(X)}
	\arrow[3]{e,t,J}{i_{k-1,k}}
\node{}
\node{}
\node{X_{k}}
	\arrow[3]{e,t,J}{i_{k}}
	\arrow{se,..}
\node{}
\node{}
\node{X}
	\arrow{se,t}{\sigma(X)}
\\
\node{}
\node{P^{m}(\G{X_{k-1}})}
	\arrow[3]{s,l}{e^{X_{k-1}}_{m}}
	\arrow[3]{e,b,J}{P^{m}(\G{i_{k-1,k}})}
\node{}
\node{}
\node{P^{m}(\G{X_{k}})}
	\arrow[3]{s,l}{e^{X_{k}}_{m}}
	\arrow[3]{e,b,J}{P^{m}(\G{i_{k}})}
\node{}
\node{}
\node{P^{m}(\G{X})}
	\arrow[3]{s,r}{e^{X}_{m}}
\\
\\
\\
\node{}
\node{X_{k-1}}
	\arrow[3]{e,t,J}{i_{k-1,k}}
\node{}
\node{}
\node{X_{k}}
	\arrow[3]{e,t,J}{i_{k}}
\node{}
\node{}
\node{X,}
\end{diagram}
\end{equation}
where $i_{k-1,k} : X_{k-1} \hookrightarrow X_{k}$ denotes the canonical inclusion.
Also the $k'$-invariant $f_{k-1}$ induces the following cofibration sequence:
$$
{S_{k-1}} \overset{f_{k-1}}\to {X_{k-1}} \overset{i_{k-1,k}}\hookrightarrow {X_{k}} \overset{g_{k}}\to \Sigma{S_{k-1}}.
$$
%
The obstruction to extend $\sigma_{k-1}(X)$ to $X_{k}$ is given by a map $\gamma = P^m(\G{i_{k-1,k}}){\comp}\sigma_{k-1}(X){\comp}f_{k-1} : S_{k-1} \to P^{m}(\G{X_{k}})$.
But the commutativity of the diagram (\ref{eq:cat-structure}) implies that $p^{\G{X_{k}}}_m{\comp}\gamma = i_{k-1,k}{\comp}e^{X_{k-1}}_m{\comp}\sigma_{k-1}(X){\comp}f_{k-1} \sim i_{k-1,k}{\comp}f_{k-1} \sim 0$.
Hence $\gamma$ has a unique lift $\hat\gamma : S_{k-1} \to E^{m+1}(\G{X_{k}})$, which vanishes in $E^{m+1}(\G{X})$, since $P^m(\G{i_{k}}){\comp}\gamma = P^m(\G{i_{k-1}}){\comp}\sigma_{k-1}(X){\comp}f_{k-1} \sim \sigma(X){\comp}i_{k-1}{\comp}f_{k-1} \sim 0$.
Since $X$ and $X_{k}$ are simply connected and $(X,X_{k})$ is $k$-connected, $(E^{m+1}(\G{X}),E^{m+1}(\G{X_{k}}))$ is ($k + 2m -1$)-connected.
Hence $\gamma$ vanishes since the dimension of $S_{k-1}$ is at most $k \leq k+2m-2$.
Thus there is a map $\sigma'_{k} : X_{k} \to P^{m}(\G{X_{k}})$ such that $\sigma'_{k}{\comp}i_{k-1,k} \sim P^m(\G{i_{k-1,k}}){\comp}\sigma_{k-1}(X)$.

The difference between $e^{X_{k}}_m{\comp}\sigma'_{k}$ and the identity $1_{X_{k}}$ with respect to the co-action of ${\Sigma}S_{k-1}$ is given by a map $\delta : {\Sigma}S_{k-1} \to X_{k}$.
By (\ref{eq:Sugawara-exact}), with $V=S_{k-1}$, $X=X_{k-1}$ and $W=X_{k}$, we have a surjection $
[{\Sigma}S_{k-1},{\Sigma}\G{X_{k}}] \to [{\Sigma}S_{k-1},X_{k}]
$ and hence $\delta$ can be pulled back to $\delta_0 : \Sigma{S_{k-1}} \to {\Sigma}\G{X_{k}} \subset P^m(\G{X_{k}}) $.
By adding $\delta_0$ to $\sigma'_k$, we have a compression $\sigma_k$ of $1_{X_{k}}$ which is an extension of $\sigma_{k-1}(X)$.
Thus $\cat{X_{k}} \leq m$.

The difference between $P^{m}(\G{i_k}){\comp}\sigma_k$ and $\sigma(X){\comp}i_{k}$ with respect to the co-action of $\Sigma{S_{k-1}}$ is given by a map $\varepsilon : {\Sigma}S_{k-1} \to P^{m}(\G{X})$ which vanishes in $P^{\infty}(\G{X}) \simeq X$.
Thus $\varepsilon$ can be lifted to $\varepsilon_0 : {\Sigma}S_{k-1} \to E^{m+1}(\G{X})$.
Since $(E^{m+1}(\G{X}),E^{m+1}(\G{X_{k}}))$ is ($k + 2m -1$)-connected and the dimension of ${\Sigma}S_{k-1}$ is at most $k+1 \leq k+2m-1$, $\varepsilon_0$ can be compressed into $E^{m+1}(\G{X_{k}})$; $\varepsilon_0 : {\Sigma}S_{k-1} \overset{\varepsilon'_0}\to E^{m+1}(\G{X_k}) \subset E^{m+1}(\G{X})$.
Again by adding $p^{\G{X_{k}}}_{m}{\comp}\varepsilon'_0$ to $\sigma_k$, we obtain a new structure map $\sigma_k(X) : X_{k} \to P^{m}(\G{X_{k}})$ for $\cat{X_{k}} \leq m$ as a compression of $\sigma(X){\comp}{i_k}$, which is also an extension of $\sigma_{k-1}(X)$.
Thus we obtain the compression $\sigma_k(X)$ of $\sigma(X){\comp}{i_k}$ for all $k$.
The latter part of the theorem is clear by the definition of $\sigma_{k}(X)$.
\end{Proof}

Now we apply the results in Section \ref{sect:hopf-cat} for homology decompositions.
Then we can show Theorems \ref{thm:stably-trivial} (1) and \ref{thm:stably-non-trivial} (1) in a slightly stronger form.
Let $\cat{X_{k}} = m \geq 1$ for some $k \geq 1$, and let $f_k(X) : S_k(X) \to X_{k}$ be the $k'$-invariant of the $k$-th stage.
We show that the obstruction for $X_{k+1}$ to satisfy $\cat{X_{k+1}} = m$ is the set of higher Hopf invariants $H^S_m(f_k(X))$.
\begin{Thm}\label{thm:obstruction-beta}
$\cat{X_{k+1}} = m$ if and only if $E^{m+1}(\G{i_{k,k+1}})_{\ast}H^S_m(f_k(X)) \ni 0$.
Moreover, if one of the following three conditions is satisfied, then $\cat{X_{k+1}} = m$ if and only if $H^S_m(f_k(X)) \ni 0$.
\begin{description}
\item{(i)} $m\geq 3$.
\item{(ii)} $X$ is simply connected and $m\geq 2$.
\item{(iii)} $X$ is simply connected and $\Ext(H_{k+1}(X),H_2(X){\otimes}H_{k+1}(X))$ $=$ $0$.
\end{description}
\end{Thm}
\begin{Proof}
If $\cat{X_{k+1}} = m$, then by Theorem \ref{thm:hom-decomp}, there exists a structure map $\sigma_{k}(X)$ for $\cat{X_{k}} = m$ such that $P^{m}(\G{i_{k,k+1}}){\comp}\sigma_{k}(X)$ is extendible to $X_{k+1}$.
For this particular choice of $\sigma_k(X)$, we obtain, by Proposition \ref{prop:obstruction}, the following diagram except for the dotted arrows commutative up to homotopy.
\begin{equation*}
\begin{diagram}
\dgHORIZPAD=.5em
\dgVERTPAD=1ex
\node{}
\node{S_k(X)}
        \arrow{e,t}{f_k(X)}
        \arrow{s,l}{H^{\sigma_k(X)}_{m}(f_k(X))}
\node{X_k}
        \arrow{e,t,J}{i_{k,k+1}}
        \arrow{s,r}{\sigma_k(X)}
\node{X_{k+1}}
        \arrow[2]{s,r,..}{\lower-6.5ex\hbox{${}_{\sigma_{k+1}}$}}
\\
\node{\G{X_k}}
        \arrow{e,J}
        \arrow{s,l,L}{\G{i_{k,k+1}}}
\node{E^{m+1}(\G{X_k})}
        \arrow{e,t,..}{p^{\G{X_{k}}}_{m}}
        \arrow{s,l,L}{E^{m+1}(\G{i_{k,k+1}})}
\node{P^{m}(\G{X_k})}
        \arrow{s,r,L}{P^{m}(\G{i_{k,k+1}})}
\node{}
\\
\node{\G{X_{k+1}}}
        \arrow{e,J}
\node{E^{m+1}(\G{X_{k+1}})}
        \arrow{e,b}{p^{\G{X_{k+1}}}_{m}}
\node{P^{m}(\G{X_{k+1}})}
        \arrow{e,b,J}{\iota^{\G{X_{k+1}}}_{m}}
\node{P^{m+1}(\G{X_{k+1}}).}
\end{diagram}
\end{equation*}
Thus the extendibility of $P^{m}(\G{i_{k,k+1}}){\comp}\sigma_{k}(X)$ implies that $p^{\G{X_{k+1}}}_{m}{\comp}E^{m+1}(\G{i_{k,k+1}}){\comp}H^{\sigma_k(X)}_{m}(f_k(X))$ $\sim$ $P^{m}(\G{i_{k,k+1}}){\comp}\sigma_{k}(X){\comp}f_k(X)$ is trivial.
Since $p^{\G{X_{k+1}}}_{m}$ induces a split monomorphism of homotopy groups with any coefficient groups, we have that $E^{m+1}(\G{i_{k,k+1}}){\comp}H^{\sigma_k(X)}_{m}(f_k(X))$ is trivial.
In each case of (i), (ii) or (iii), we show $\pi_{k}(E^{m+1}(\G{X_{k+1}}),E^{m+1}(\G{X_k});H_{k+1}(X))$ $=$ $\pi_{k+1}(E^{m+1}(\G{X_{k+1}}),E^{m+1}(\G{X_k});H_{k+1}(X))$ $=$ $0$.

Case (i):
The pair ($E^{m+1}(\G{X_{k+1}}),E^{m+1}(\G{X_k})$) is ($k+m-1$)-connected, $m \geq 3$.
Hence, $\pi_{k+1}(E^{m+1}(\G{X_{k+1}}),E^{m+1}(\G{X_k});H_{k+1}(X)) = 0$ for dimensional reasons.

Cases (ii) and (iii):
The pair ($E^{m+1}(\G{X_{k+1}}),E^{m+1}(\G{X_k})$) is ($k+2m-1$)-connected.
When $m \geq 2$, we obtain $\pi_{k+1}(E^{m+1}(\G{X_{k+1}}),E^{m+1}(\G{X_k});H_{k+1}(X)) = 0$ for dimensional reasons.
When $m=1$, by the Universal Coefficient Theorem for homotopy groups, we obtain 
\begin{align*}&
\pi_{k+1}(E^{2}(\G{X_{k+1}}),E^{2}(\G{X_k});H_{k+1}(X)) 
= \Ext(H_{k+1}(X),H_{k+2}(E^{2}(\G{X_{k+1}}),E^{2}(\G{X_k}))).
\end{align*}
For dimensional reasons, we have 
\begin{align*}&
H_{k+2}(E^{2}(\G{X_{k+1}}),E^{2}(\G{X_k})) \cong H_{k+2}(S_1(X){\ast}S_k(X){\vee}S_k(X){\ast}S_1(X)) 
\\&\qquad
\cong H_2(X){\otimes}H_{k+1}(X){\oplus}H_{k+1}(X){\otimes}H_2(X),
\quad\text{and hence}\\&
\pi_{k+1}(E^{2}(\G{X_{k+1}}),E^{2}(\G{X_k});H_{k+1}(X)) 
\\&\qquad
\cong \Ext(H_{k+1}(X),H_2(X){\otimes}H_{k+1}(X)){\oplus}\Ext(H_{k+1}(X),H_{k+1}(X){\otimes}H_2(X)),
\end{align*}
which is trivial if $\Ext(H_{k+1}(X),H_2(X){\otimes}H_{k+1}(X))=0$.

Hence by assuming (i), (ii) or (iii), we obtain that $E^{m+1}(\G{i_{k,k+1}})_{\ast}$ :  $\pi_{k}(E^{m+1}(\G{X_{k}}))$ = $[S_k(X),E^{m+1}(\G{X_k})]$ $\to$ $[S_k(X),E^{m+1}(\G{X_{k+1}})]$ = $\pi_{k}(E^{m+1}(\G{X_{k+1}}))$ has no non-trivial kernel.
Thus the set $H^S_m(f_k(X))$ contains $0=H^{\sigma_k(X)}_{m}(f_k(X))$.
The converse is an immediate consequence of Theorem \ref{thm:stably-trivial} (1).
\end{Proof}

Let $\cat{X_{k+1}} = m+1$, in other words, the set $E^{m+1}(\G{i_{k,k+1}})_{\ast}H^S_m(f_k(X))$ does not contain $0$.
Then the following theorem is an immediate consequence of Theorem \ref{thm:stably-trivial} (2).
\begin{Thm}\label{thm:obstruction-beta'}
If the set ${\Sigma}^n_{\ast}H^S_m(f_k(X))$ contains $0$, then $\cat{X_{k+1}{\times}S^{n}} = \cat{X_{k+1}} = m+1$.
\end{Thm}
\begin{Cor}\label{cor:obstruction-beta'}
If the set ${\mathcal H}^S_m(f_k(X))$ contains $0$, then $X_{k+1}$ is a counter example to Ganea's conjecture.
\end{Cor}

Also let $n \geq 1$, $\cat{Y} \geq n+1$.
Let $\cat{Y_{h}} = n \geq 1$ for some $h \geq 1$, and let $f_h(Y) : S_h(Y) \to Y_{h}$ be the $k'$-invariant of the $h$-th stage.
We know, by Theorem \ref{thm:obstruction-beta}, that the obstruction for $Y_{h+1}$ to satisfy $\cat{Y_{h+1}} \leq n$ is the set of Hopf invariants $H^{S}_{n}(f_h(Y))$.
We define another set 
\begin{Def}\label{def:set-higher-hopf2}\quad
\begin{equation*}
H^{S}_{m,n}(f_k(X),f_h(Y)) = \{g_X{\ast}g_Y \,\vert\, g_X \in H^{S}_m(f_k(X)) ~\text{and}~ g_Y \in H^{S}_n(f_h(Y))\}
\end{equation*}
\end{Def}
Then we have the following theorem.
\begin{Thm}\label{thm:obstruction-beta'2}
If the set $H^{S}_{m,n}(f_k(X),f_h(Y))$ contains $0$, then $\cat{X_{k+1}{\times}Y_{h+1}} < m+n+2$ $=$ $\cat{X_{k+1}} + \cat{Y_{h+1}}$.
\end{Thm}
\begin{Proof}
The proof is obtained by a similar argument given in the proof of Theorem \ref{thm:stably-trivial} (2) using the following diagram instead of the diagram in Proposition \ref{prop:obstruction2}:

\begin{equation*}
\divide\dgARROWLENGTH by2
\dgHORIZPAD=.5em
\dgVERTPAD=1ex
\begin{diagram}
\node{\substack{S_{k+h+1}(X_{k+1}{\times}Y_{h+1})}}
        \arrow{s,r}{H^{\sigma_k(X)}_{m+1}(f_k(X)){\ast}H^{\sigma_h(Y)}_{n+1}(f_h(Y))}
        \arrow[4]{e,t}{[f_k(X),f_h(Y)]^r}
\node{}
\node{}
\node{}
\node{\substack{X_{k+1}{\times}Y_{h}{\cup}X_{k}{\times}Y_{h+1}}}
        \arrow{s,r}{(\sigma_{k+1}(X){\times}\sigma_{h+1}(Y)){\vert}}
        \arrow[3]{e,J}
\node{}
\node{}
\node{\substack{X_{k+1}{\times}Y_{h+1}}}
        \arrow[2]{s,..}
\\
\node{\substack{E^{m+1}(\G{X_k}){\ast}E^{n+1}(\G{Y_h})}}
        \arrow{s,r,J}{E^{m+1}(\G{i_{k,k+1}}){\ast}E^{n+1}(\G{i_{h,h+1}})}
        \arrow[4]{e,t,..}{[p^{\G{X_k}}_m,p^{\G{Y_h}}_n]^r}
\node{}
\node{}
\node{}
\node{\substack{P^{m+1}(\G{X_{k+1}}){\times}P^{n}(\G{Y_{h}})\\~\qquad{\cup}P^{m}(\G{X_{k}}){\times}P^{n+1}(\G{Y_{h+1}})}}
        \arrow{s,J}
\\
\node{\substack{E^{m+1}(\G{X_{k+1}}){\ast}E^{n+1}(\G{Y_{h+1}})}}
        \arrow[4]{e,t}{[p^{\G{X_{k+1}}}_m,p^{\G{Y_{h+1}}}_n]^r}
\node{}
\node{}
\node{}
\node{\substack{P^{m+1}(\G{X_{k+1}}){\times}P^{n}(\G{Y_{h+1}})\\~\qquad{\cup}P^{m}(\G{X_{k+1}}){\times}P^{n+1}(\G{Y_{h+1}})}}
        \arrow[3]{e,J}
\node{}
\node{}
\node{\substack{P^{m+1}({\G{X_{k+1}}}){\times}P^{n+1}({\G{Y_{h+1}}}),}}
\end{diagram}
\end{equation*}
The details are left to the reader.
\end{Proof}

\section{Higher Hopf invariants for some examples}\label{sect:hhi-example}

In this section, we compute the higher Hopf invariants for well-known examples, which yields a generalisation of the main result of \cite{Iwase:counter-ls}:
We denote by $\complex \homeo \real^2$ the field of complex numbers, by $\quaternion \homeo \real^4$ the algebra of quaternion numbers and by $\cayley \homeo \real^8$ the Cayley algebra:

\begin{Exam}\label{exam:complex}
We know that $\cat{{\complex}P^m} = m$ and $\dim{{\complex}P^m} = 2m$.
Hence $\dim{{\complex}P^m} = 2m \leq {2}m + 2 - 2$.
Thus ${\complex}P^m$ has a unique structure for $\cat{{\complex}P^m} = m$.
The higher Hopf invariant $H_m$ $:$ $[V,{\complex}P^m] \to [V,\Omega{\complex}P^m{\ast}{\cdots}{\ast}\Omega{\complex}P^m]$ gives a (unique) homomorphism
\begin{equation*}
H_m : \pi_{2m+1}({\complex}P^m) \to \pi_{2m+1}(\Omega{\complex}P^m{\ast}{\cdots}{\ast}\Omega{\complex}P^m) \cong \integer
\end{equation*}
with the canonical projection $p^{S^1}_m$ $:$ $S^{2m+1}$ $\to$ ${\complex}P^m$ a `higher Hopf invariant one' element.
\end{Exam}

\begin{Exam}\label{exam:quaternion}
We know that $\cat{{\quaternion}P^m} = m$ and $\dim{{\quaternion}P^m} = 4m$.
Hence $\dim{{\quaternion}P^m} = 4m \leq {4}m + 4 - 2$.
Thus ${\quaternion}P^m$ has a unique structure for $\cat{{\quaternion}P^m} = m$.
The higher Hopf invariant $H_m$ $:$ $[V,{\quaternion}P^m] \to [V,\Omega{\quaternion}P^m{\ast}{\cdots}{\ast}\Omega{\quaternion}P^m]$ gives a (unique) homomorphism
\begin{equation*}
H_m : \pi_{4m+3}({\quaternion}P^m) \to \pi_{4m+3}(\Omega{\quaternion}P^m{\ast}{\cdots}{\ast}\Omega{\quaternion}P^m) \cong \integer
\end{equation*}
with the canonical projection $p^{S^3}_m$ $:$ $S^{4m+3}$ $\to$ ${\quaternion}P^m$ a `higher Hopf invariant one' element.
\end{Exam}

\begin{Exam}\label{exam:cayley}
We know that $\cat{{\cayley}P^2} = 2$ and $\dim{{\cayley}P^2} = 16$.
Hence $\dim{{\cayley}P^2} = 16 \leq {8}\times{2} + 8 - 2$.
Thus ${\cayley}P^2$ has a unique structure for $\cat{{\cayley}P^2} = 2$.
The higher Hopf invariant $H_2$ $:$ $[V,{\cayley}P^2] \to [V,\Omega{\cayley}P^2{\ast}\Omega{\cayley}P^2{\ast}\Omega{\cayley}P^2]$ gives a (unique) homomorphism
\begin{equation*}
H_2 : \pi_{23}({\cayley}P^2) \to \pi_{23}(\Omega{\cayley}P^2{\ast}\Omega{\cayley}P^2{\ast}\Omega{\cayley}P^2) \cong \integer.
\end{equation*}
But there are no elements of `higher Hopf invariant one':
The existence of such a higher Hopf invariant one element implies that the Hopf space $S^7$ is homotopy associative.
As is well-known, the $p$-local Hopf space $S^7_{(p)}$ is homotopy associative for $p \geq 5$ (in view of \cite{Sugawara:group-like}, \cite{Stasheff:higher-associativity} and \cite{IM:higher-associativity}, it actually is an $A_{p-1}$-space).
However, by using primary cohomology operations, one can easily see that any Hopf structure on $S^7_{(3)}$ is not homotopy associative.
Hopf space theorists were, however, much more interested in the case $p=2$.
And it was known by Goncalves \cite{Goncalves:non-hom-ass-S7}, using higher order cohomology operations, and by Hubbuck \cite{Hubbuck:two-example}, using K-theory Adams operations, that any Hopf structure on $S^7_{(2)}$ is not homotopy associative (but the result itself had already been known by James).
Hence, the image of the higher Hopf invariant homomorphism is in $6\integer \subset \integer$.
\end{Exam}

\begin{Exam}\label{exam:general-counter-ls}
For $m,n \geq 1$ and $p \geq m+2$, let $f_{m,p} = p^{S^1}_{m}{\comp}g_{m,p}$ $:$ $S^{2mp+2(p-1)-1}$ $\to$ ${\complex}P^m$, where $g_{m,p}$ $:$ $S^{2mp+2(p-1)-1}$ $\to$ $S^{2m+1}$ denotes the generator of $\pi_{2mp+2(p-1)-1}(S^{2m+1}) \cong \integer/p\integer$ and $p^{S^1}_m$ $:$ $S^{2m+1}$ $\to$ ${\complex}P^m$ denotes the projection which gives a `higher Hopf invariant one' element as in Example \ref{exam:complex}.
For dimensional reasons, the map $g_{m,p}$ is a co-H-map.
Hence, by Proposition \ref{prop:fundamental-property}(1), $H_m(f_{m,p}) = g_{m,p} \ne 0$, and hence ${\Sigma}^{n}H_m(f_{m,p})$ is trivial if and only if $n \geq 2$, by Theorem 13.4 in Toda \cite{Toda:composition-methods}.

Let $Q_{m,p}$ be the mapping cone of $f_{m,p}$: $Q_{m,p} = {\complex}P^m \cup_{f_{m,p}} e^{2mp+2(p-1)}$.
Then, by Theorems \ref{thm:stably-trivial} and \ref{thm:stably-non-trivial} and Remark \ref{rem:remove-inc}, it follows that $\cat{Q_{m,p}} = m+1$ and $m+1 \leq \cat(Q_{m,p}{\times}S^1) \leq m+2$ but $\cat(Q_{m,p}{\times}S^n) = m+1$ for $n \geq 2$ 
by Theorem \ref{thm:stably-trivial}(2).
\end{Exam}
\begin{Rems}
\begin{description}
\item{(1)}
Every example in this section supports Conjecture \ref{conj:iwase-number}.
\item{(2)}
The space $Q_{m,p}$ in Example \ref{exam:general-counter-ls} is a generalisation of $Q_p$ in \cite{Iwase:counter-ls} except $Q_2$.
Actually, Ganea's conjecture for $Q_{m,p}$ is true if we consider $\catq{}$, for $q \ne p$, instead of $\cat{}$ or $\catp{}$.
\end{description}
\end{Rems}

\section{LS category of sphere-bundles-over-spheres}\label{sect:cat-sphere-bundle-over-sphere}

Let $r \geq 1$, $t \geq 1$ and $E$ be a fibre bundle over $S^{t+1}$ with fibre $S^{r}$.
Then $E$ has a CW decomposition $S^{r} \cup_{\alpha} e^{t+1} \cup_{\psi} e^{t+r+1}$ with $\alpha : S^{t} \to S^{r}$ and $\psi : S^{t+r} \to Q$, $Q = S^{r} \cup_{\alpha} e^{t+1}$.
We identify $H^S_1(\alpha)$ with its unique element, say $H_1(\alpha)$, since $S^k$ has a unique structure map for $\cat{S^k} = 1$.

\begin{Fact}\label{fact:cat-Q-0}
Let $\alpha = {\pm}1_{S^{r}}$, the identity.
Then clearly $\cat{Q} = 0$ and $\cat{E} = 1$.
In addition, $\cat{Q{\times}S^n} = 1$ and $\cat{E{\times}S^n} = 2$ for $n \geq 1$.
\end{Fact}

The following fact is an immediate consequence of Berstein-Hilton \cite{BH:category} and a cup length consideration.
\begin{Fact}\label{fact:cat-Q-1}
Let $\alpha \ne {\pm}1_{S^{r}}$.
Hence $1 \leq \cat{Q} \leq 2$.
Then $\cat{Q} = 2$ if and only if $H_1(\alpha) \ne 0$.
In particular if $H_1(\alpha) = 0$, we can easily obtain that $\cat{Q} = 1$ and $\cat{E} = 2$.
In this case, it also follows that $\cat{Q{\times}S^n} = 2$ and $\cat{E{\times}S^n} = 3$ for $n \geq 1$.
\end{Fact}

By Theorem \ref{thm:stably-non-trivial} and Remark \ref{rem:remove-inc}, we can extend the main result of \cite{Iwase:counter-ls}.

\begin{Thm}\label{thm:cat-Q-2}
Let $H_1(\alpha) \ne 0$.
Hence $\cat{Q} = 2$.
Then for $n \geq 1$, $\cat{Q{\times}S^n} = 3$ if ${\Sigma}^{n}H_1(\alpha) \ne 0$ with $n \geq t-2r+2$ or ${\Sigma}^{n+1}H_1(\alpha) \ne 0$.
\end{Thm}

We give a partial answer to Ganea's Problem 4 (see \cite{Ganea:conjecture}) for sphere-bundles-over-spheres.
To show this, we need the following lemma.

\begin{Lem}\label{lem:trivial}
The collapsing map $q : E \to E/Q = S^{t+r+1}$ induces a map with trivial kernel
\begin{equation*}
({\Sigma}^nq)^{\ast} : [S^{n+r+t+1},S^{n+kr}] \to [{\Sigma}^nE,S^{n+kr}]\quad\text{for all $k \geq 3$ and $n \geq 0$.}
\end{equation*}
\end{Lem}
\begin{Proof}
The cofibration sequences $S^{r+t} \overset{\psi}\to Q \overset{j}\to E$ and $S^{t} \overset{\alpha}\to S^{r} \overset{i}\to Q$ together with the bundle projection $p : E \to S^{t+1}$ induce the following commutative diagram:

\begin{equation*}
\dgHORIZPAD=.5em
\dgVERTPAD=1ex
\begin{diagram}
\node{}
\node{[S^{n+t+2},S^{n+kr}]}
        \arrow{sw,l}{{\Sigma}^{n+1}p^{\ast}}
        \arrow{s,r}{({\Sigma}^{n+1}p\vert_{{\Sigma}^{n+1}Q})^{\ast}}
\\
\node{[{\Sigma}^{n+1}E,S^{n+kr}]}
        \arrow{e,b}{{\Sigma}^{n+1}j^{\ast}}
\node{[{\Sigma}^{n+1}Q,S^{n+kr}]}
        \arrow{e,b}{{\Sigma}^{n+1}\psi^{\ast}}
        \arrow{s,r}{{\Sigma}^{n+1}i^{\ast}}
\node{[S^{n+r+t+1},S^{n+kr}]}
        \arrow{e,b}{({\Sigma}^nq)^{\ast}}
\node{[{\Sigma}^nE,S^{n+kr}]}
\\
\node{}
\node{[S^{n+r+1},S^{n+kr}],}
\end{diagram}
\end{equation*}
where the column and row are exact sequences.
Since $k \geq 3$, we have $n+r+1<n+kr$ and $\pi_{n+r+1}(S^{n+kr}) = 0$.
Hence $({\Sigma}^{n+1}p\vert_{{\Sigma}^{n+1}Q})^{\ast}$ is surjective, and so is $({\Sigma}^{n+1}j)^{\ast}$.
Thus $({\Sigma}^{n+1}\psi)^{\ast}$ is trivial and the map $({\Sigma}^nq)^{\ast}$ has trivial kernel for $k \geq 3$ and $n \geq 0$.
\end{Proof}

\begin{Thm}\label{thm:cat-E-2}
Let $H_1(\alpha) \ne 0$.
Hence $2 \leq \cat{E} \leq 3$.
Then $\cat{E} = 3$ if ${\Sigma}^{r+1}h_2(\alpha) \ne 0$.
Also $\cat{E} = 2$ if $H^S_2(\psi)$ or its subset $H^{SS}_2(\psi)$ (see Remark \ref{rem:category-W} for its definition) contains $0$.
\end{Thm}
\begin{Rem}
In the latter case of $\cat{E} = 2$, it is known that $\cat{E{\times}S^n} = 3$ for $n \geq 1$ by using a cup length argument on the cohomology ring (see Singhof \cite{Singhof:minimal-cover}).
\end{Rem}
\begin{Proof*}{\it Proof of Theorem \ref{thm:cat-E-2}.}
Let $q' : E \to E/S^{r}$, $q'' : Q \to Q/S^{r}$ be respectively the collapsing maps.
The reduced diagonal map $\hat{\Delta}_2 : E \to E{\wedge}E$ factors as $E \overset{q'}\to S^{t+1} {\cup} e^{r+t+1} \overset{\hat{\Delta}_2}\to (S^{r} \cup e^{t+1}){\wedge}(S^{r} \cup e^{t+1}) \subset E{\wedge}E$, which is an extension of the map $Q \overset{q''}\to S^{t+1} \overset{\Sigma{h_2(\alpha)}}\to S^{r}{\wedge}S^{r} \subset E{\wedge}E$ by Theorem 5.14 of Boardman-Steer \cite{BS:hopf-invariants}.
In this case, the generator in $H^{r+t+1}(E;\integer)$ is a cup product of generators in $H^{r}(E;\integer)$ and $H^{t+1}(E;\integer)$.
Then it follows that the mapping degrees of $\hat{\Delta}_2$ on $e^{r}{\wedge}e^{t+1}$ and $e^{t+1}{\wedge}e^{r}$ are 1.
Hence the reduced diagonal map $\hat{\Delta}_3 = (\hat{\Delta}_2{\wedge}1_E){\comp}\hat{\Delta}_2 : E \to E{\wedge}E{\wedge}E$ factors as $E \overset{q}\to S^{r+t+1} \overset{\simeq}\to S^{t+1}{\wedge}S^{r} \overset{{\Sigma}^{r+1}h_2(\alpha)}\to S^{r}{\wedge}S^{r}{\wedge}S^{r} \subset E{\wedge}E{\wedge}E$.
Hence $\hat{\Delta}_3$ is the composition of ${\Sigma}^{r+1}h_2(\alpha){\comp}q$ with a suitable inclusion, which does not depend on the choice of $\sigma(Q)$.
By Lemma \ref{lem:trivial}, $q^{\ast} : [S^{r+t+1},S^{3r}] \to [E,S^{3r}]$ has trivial kernel, and hence the non-triviality of ${\Sigma}^{r+1}h_2(\alpha)$ implies the non-triviality of ${\Sigma}^{r+1}h_2(\alpha){\comp}q$.
Then, for dimensional reasons, it follows that $\hat{\Delta}_3$ is also non-trivial.
Then by Theorem \ref{thm:reduced-diagonal}, it follows that $H^{S}_{2}(\psi)$ does not contain $0$, 
and hence we see that $\cat{E} = 3$ by Theorem \ref{thm:obstruction-beta}.
The latter part is clear by Theorem \ref{thm:obstruction-beta}.
\end{Proof*}

We next study the LS category of $E{\times}S^n$.
To do this, we need the following lemma.

\begin{Lem}\label{lem:trivial'}
The collapsing map $\hat{q} : E{\times}S^n \to E{\times}S^n/(E \cup Q{\times}S^n) = S^{n+r+t+1}$ induces a map with trivial kernel
\begin{equation*}
\hat{g}^{\ast} : [S^{n+r+t+1},S^{n+kr}] \to [E{\times}S^n,S^{n+kr}]\quad\text{for all $k \geq 3$ and $n \geq 1$.}
\end{equation*}
\end{Lem}
\begin{Proof}
Let us recall that the space $E \vee S^n$ is a retractile subspace (see Zabrodsky \cite{Zabrodsky:Hopf-space}) of both $E{\times}S^n$ and $E \cup Q{\times}S^n$.
The cofibration sequences ${S^{n+r+t}} \overset{\psi'}\to {E \cup Q{\times}S^n} \overset{j'}\hookrightarrow {E{\times}S^n}$, ${E \vee S^n} \overset{i_0}\hookrightarrow {E{\times}S^n} \overset{q_0}\to {E{\wedge}S^n}$ and ${E \vee S^n} \overset{i'_0}\hookrightarrow {E \cup Q{\times}S^n} \overset{q'_0}\to {Q{\wedge}S^n}$ induce the following commutative diagram:
\begin{equation*}
\divide\dgARROWLENGTH by2
\dgHORIZPAD=.5em
\dgVERTPAD=1ex
\begin{diagram}
\node{[{\Sigma}(E{\wedge}S^n),S^{n+kr}]}
        \arrow{e,t}{({\Sigma}^{n+1}j)^{\ast}}
        \arrow{s,l}{({\Sigma}q_0)^{\ast}}
\node{[{\Sigma}(Q{\wedge}S^n),S^{n+kr}]}
        \arrow{e,t}{({\Sigma}^{n+1}\psi)^{\ast}}
        \arrow{s,l}{({\Sigma}q'_0)^{\ast}}
\node{[S^{n+r+t+1},S^{n+kr}]}
        \arrow{s,=}
        \arrow{e,t}{({\Sigma}^nq)^{\ast}}
\node{[E{\wedge}S^n,S^{n+kr}]}
\\
\node{[{\Sigma}(E{\times}S^n),S^{n+kr}]}
        \arrow{e,t}{({\Sigma}j')^{\ast}}
        \arrow{s,r}{({\Sigma}i_0)^{\ast}}
\node{[{\Sigma}(E \cup Q{\times}S^n),S^{n+kr}]}
        \arrow{e,t}{({\Sigma}\psi')^{\ast}}
        \arrow{s,r}{({\Sigma}i'_0)^{\ast}}
\node{[S^{n+r+t+1},S^{n+kr}]}
        \arrow{e,t}{\hat{q}^{\ast}}
\node{[E{\times}S^n,S^{n+kr}]}
\\
\node{[{\Sigma}E{\vee}{\Sigma}S^{n},S^{n+kr}]}
        \arrow{e,=}
\node{[{\Sigma}E{\vee}{\Sigma}S^{n},S^{n+kr}],}
\end{diagram}
\end{equation*}
where the columns are exact sequences and the rows are split short exact sequences with natural splittings (see Zabrodsky \cite{Zabrodsky:Hopf-space}).
By the proof of Lemma \ref{lem:trivial}, $({\Sigma}^{n+1}j)^{\ast}$ is surjective, and hence so is $({\Sigma}j')^{\ast}$.
Thus $({\Sigma}\psi')^{\ast}$ is trivial, and hence $\hat{q}^{\ast}$ has trivial kernel for $k \geq 3$ and $n \geq 0$.
\end{Proof}

\begin{Thm}\label{thm:cat-E-3}
Let ${\Sigma}^{r+1}h_2(\alpha) \ne 0$.
Hence $\cat{E} = 3$.
Then for $n \geq 1$, $\cat{E{\times}S^n} = 4$ if ${\Sigma}^{n+r+1}h_2(\alpha) \ne 0$.
Also for $n \geq 1$, $\cat{E{\times}S^n} = 3$ if $H^S_3(\psi')$, ${\Sigma}^{n}_{\ast}H^S_2(\psi)$ or ${\Sigma}^{n}_{\ast}H^{SS}_2(\psi)$ $\ni$ $0$.
\end{Thm}
\begin{Proof}
By the proof of Theorem \ref{thm:cat-E-2}, the $n$-fold suspension of the reduced diagonal map $\hat{\Delta}_3{\wedge}1_{S^n}$ $:$ $E{\wedge}S^n$ $\to$ $E{\wedge}E{\wedge}E{\wedge}S^n$ is the composition of ${\Sigma}^{n+r+1}h_2(\alpha){\comp}\hat{q}$ with a suitable inclusion, which does not depend on the choice of $\sigma(E{\cup}Q{\times}S^n)$.
By Lemma \ref{lem:trivial'}, $\hat{q}^{\ast} : [S^{n+r+t+1},S^{n+3r}] \to [E{\times}S^n,S^{n+3r}]$ has trivial kernel, and hence the non-triviality of ${\Sigma}^{n+r+1}h_2(\alpha)$ implies the non-triviality of ${\Sigma}^{n+r+1}h_2(\alpha){\comp}\hat{q}$.
Then, for dimensional reasons, it follows that $\hat{\Delta}_3{\wedge}1_{S^n}$ is also non-trivial, and hence the four-fold reduced diagonal $\hat{\Delta}'_4$ of $E{\times}S^{n}$ is non-trivial.
Thus by Theorem \ref{thm:reduced-diagonal}, ${\Sigma}_{\ast}H^{S}_{3}(\psi')$ and $H^{S}_{3}(\psi')$ do not contain $0$, 
and hence we see that $\cat{E{\times}S^n} = 4$ by Theorem \ref{thm:obstruction-beta} with $m=3$.
The latter part is clear by Theorem \ref{thm:obstruction-beta} (in the case of $H^S_3(\psi')$) and Theorem \ref{thm:obstruction-beta'} (in the case of ${\Sigma}^{n}_{\ast}H^S_2(\psi)$ or ${\Sigma}^{n}_{\ast}H^{SS}_2(\psi)$).
\end{Proof}

\section{Manifold examples}\label{sect:mfd-example}

The Hopf fibration $\sigma_4 : S^7 \to S^4$ is given as a principal $Sp(1)$-bundle.
Taking orbits of the action of $U(1) \subset Sp(1)$ on $S^7$, we obtain a fibre bundle $\complex{P}^3 \to S^4$ with fibre $Sp(1)/U(1) \homeo S^2$; the structure group $Sp(1)$ acts on the fibre $S^2$ via a map, say $\mu_0 : S^3{\times}S^2 \to S^2$.
Here the CW decomposition of $\complex{P}^3$ is known as $\complex{P}^3 = \complex{P}^2 {\cup}_{p^{S^1}_2} e^{6} = S^{2} \cup_{\eta} e^{4} {\cup}_{p^{S^1}_2} e^{6}$.
Hence the attaching map $p^{S^1}_2$ of the top cell of $\complex{P}^3$ is given by the composition
\begin{equation*}
\multiply\dgARROWLENGTH by2
\begin{diagram}
\node{S^3{\ast}S^1} \arrow{e,t}{[C(1_{S^3}),\Sigma{1_{S^1}}]^r} \node{C(S^3){\cup}S^3{\times}S^2} \arrow{e,t}{\hat{\chi}_4} \node{\complex{P}^2,}
\end{diagram}
\end{equation*}
where $\hat{\chi}_4$ denotes the map defined by $\hat{\mu}\vert_{S^3{\times}S^2} = \mu_0$ and $\hat{\mu}\vert_{C(S^3)} = \chi_4$, the characteristic map of the top cell of $\complex{P}^2$.
\begin{Def}
For any map $\beta : S^t \to S^3$, we may assume that the suspension $\Sigma{\beta} : S^{t+1} \to S^4$ is a $C^{\infty}$-map by suitably deforming it up to homotopy, since $S^{t+1}$ and $S^{4}$ are closed $C^{\infty}$-manifolds.
We define $E(\beta)$ to be the total space of the $Sp(1)$-bundle $E(\beta) \to S^{t+1}$ induced by the $C^{\infty}$-map $\Sigma{\beta}$ from the $Sp(1)$-bundle ${\complex}P^3 \to S^4$.
Hence $E(\beta)$ is an orientable, closed $C^{\infty}$-manifold with CW decomposition $E(\beta) = S^2 \cup_{\eta{\comp}\beta} e^{t+1} {\cup}_{\psi(\beta)} e^{t+3}$.
\end{Def}

For a map $\beta : S^{t} \to S^3$ and a suspension map $\gamma$ : $S^{t'} \to S^{t}$ with $3 \leq t \leq t'$, 
we denote $\beta' = \beta{\comp}\gamma$, and then we have $E(\beta') = S^2 \cup_{\eta{\comp}\beta'} e^{t'+1} {\cup}_{\psi(\beta')} e^{t'+3}$.
By putting $Q(\beta) =  S^2 \cup_{\eta{\comp}\beta} e^{t+1} \subset E(\beta)$ and $Q(\beta') =  S^2 \cup_{\eta{\comp}\beta'} e^{t'+1} \subset E(\beta')$, we have the following commutative ladder of cofibration sequences:
\begin{equation}\label{cd:h_1-beta}
\dgHORIZPAD=.5em
\dgVERTPAD=1ex
\begin{diagram}
\node{S^{t'}}
	\arrow{e,t}{\eta{\comp}\beta'}
	\arrow{s,l}{\gamma}
\node{S^2}
	\arrow{e,J}
	\arrow{s,=}
\node{Q(\beta')}
	\arrow{e}
	\arrow{s,l}{\hat{\gamma}_0}
\node{S^{t'+1}}
	\arrow{s,l}{\Sigma\gamma}
	\arrow{e}
\node{S^{3},}
	\arrow{s,=}
\\
\node{S^{t}}
	\arrow{e,t}{\eta{\comp}\beta}
\node{S^2}
	\arrow{e,J}
\node{Q(\beta)}
	\arrow{e}
\node{S^{t+1}}
	\arrow{e}
\node{S^{3},}
\end{diagram}
\end{equation}
where $\hat{\gamma}_0 = \hat\gamma\vert_{Q(\beta')}$ and $\hat{\gamma} : E(\beta') \to E(\beta)$ is the bundle map induced from $\Sigma\gamma : S^{t'+1} \to S^{t+1}$.
To compare the higher Hopf invariant of $\psi(\beta)$ with that of $\psi(\beta')$, we show the following proposition.
\begin{Prop}\label{prop:primitivity}
With respect to `standard' structure maps (see Remark \ref{rem:category-W}), $\hat{\gamma}_0$ is `primitive' in the sense of Berstein and Hilton.
\end{Prop}
\begin{Rem}
${\complex}P^2$ has a unique structure map $\sigma_2({\complex}P^2)$ for $\cat{{\complex}P^2} \leq 2$ (see Example \ref{exam:complex}).
\end{Rem}
\begin{Proof*}{Proof of Proposition \ref{prop:primitivity}.}
By Proposition \ref{prop:obstruction} and Remark \ref{rem:category-W'}, there are maps $\sigma'_2(Q(\beta)) : Q(\beta) \to P^2(\G{Q(\beta)})$ and $\sigma'_2(Q(\beta')) : Q(\beta') \to P^2(\G{Q(\beta')})$.
Since $\gamma$ is a suspension map, the `naturality' of Lemma \ref{prop:obstruction-natural} implies a homotopy (relative to $S^2$):
\begin{equation}\label{eq:sigma'}
P^2(\G{\hat{\gamma}_0}){\comp}\sigma'_2(Q(\beta')) \sim \sigma'_2(Q(\beta)){\comp}\hat{\gamma}_0 : (Q(\beta'),S^2) \to (\Sigma\G{Q(\beta')},\Sigma\G{S^2}).
\end{equation}

Let us recall what is in Remark \ref{rem:category-W}:
The difference between $e^{Q(\beta)}_2{\comp}\sigma'_2(Q(\beta))$ and $1_{Q(\beta)}$ with respect to the co-action of $S^{t+1}$ is given by a map $\delta : S^{t+1} \to Q(\beta)$ which can be pulled back to a map $\delta_0 : S^{t+1} \to {\Sigma}\G{Q(\beta)}$.
Also the difference between $e^{Q(\beta')}_2{\comp}\sigma'_2(Q(\beta'))$ and $1_{Q(\beta')}$ with respect to the co-action of $S^{t'+1}$ is given by a map $\delta' : S^{t'+1} \to Q(\beta')$ which can be pulled back to a map $\delta'_0 : S^{t'+1} \to {\Sigma}\G{Q(\beta')}$.
Let $\sigma_2(Q(\beta)) = \sigma'_2(Q(\beta)) + \iota^{\G{Q(\beta)}}_1{\comp}\delta_0$ and $\sigma_2(Q(\beta')) = \sigma'_2(Q(\beta')) + \iota^{\G{Q(\beta')}}_1{\comp}\delta'_0$, where the addition is induced from the co-actions of $S^{q+1} = Q(\beta)/S^2$ on $Q(\beta)$ and of $S^{t'+1} = Q(\beta')/S^2$ on $Q(\beta')$, respectively.
Then $\sigma_2(Q(\beta))$ and $\sigma_2(Q(\beta'))$ are genuine compressions of $1_{Q(\beta)}$ and $1_{Q(\beta')}$.

Using the homotopy (\ref{eq:sigma'}), we obtain a homotopy (relative to $S^2$)
\begin{align*}&
\hat{\gamma}_0{\comp}e^{Q(\beta')}_2{\comp}\sigma'_2(Q(\beta'))
\sim 
e^{Q(\beta)}_2{\comp}P^2(\G{\hat{\gamma}_0}){\comp}\sigma'_2(Q(\beta')) 
\sim 
e^{Q(\beta)}_2{\comp}\sigma'_2(Q(\beta)){\comp}\hat{\gamma}_0
\\&\qquad\qquad
 : (Q(\beta'),S^2) \to (Q(\beta),S^2),
\\
\intertext{and hence a homotopy (relative to $S^2$)}&
e^{Q(\beta)}_2{\comp}\sigma'_2(Q(\beta)){\comp}\hat{\gamma}_0 + \delta{\comp}{\Sigma}\gamma
\sim 
(e^{Q(\beta)}_2{\comp}\sigma'_2(Q(\beta)) + \delta){\comp}\hat{\gamma}_0
\sim 
\hat{\gamma}_0
\sim 
\hat{\gamma}_0{\comp}e^{Q(\beta')}_2{\comp}\sigma_2(Q(\beta'))
\\&\qquad
\sim
\hat{\gamma}_0{\comp}e^{Q(\beta')}_2{\comp}\sigma'_2(Q(\beta')) + \hat{\gamma}_0{\comp}\delta'
\sim e^{Q(\beta)}_2{\comp}\sigma'_2(Q(\beta)){\comp}\hat{\gamma}_0 + \hat{\gamma}_0{\comp}\delta' : (Q(\beta'),S^2) \to (Q(\beta),S^2).
\end{align*}
Hence the difference of $\delta{\comp}{\Sigma}\gamma$ and $\hat{\gamma}_0{\comp}\delta'$ is trivial in $\pi_{t'+1}(Q(\beta))$ by using the ordinary obstruction theory (see \cite{Whitehead:elements}).
Thus the difference of $\delta_0{\comp}{\Sigma}\gamma$ and ${\Sigma}\G{\hat{\gamma}_0}{\comp}\delta'_0$ in $\pi_{t'+1}({\Sigma}\G{Q(\beta)})$, which is given by a map $\epsilon_0 : S^{t'+1} \to {\Sigma}\G{Q(\beta)}$, vanishes in $\pi_{t'+1}(Q(\beta'))$.
Thus $\epsilon_0$ can be lifted uniquely to a map $\hat{\epsilon}_0 \in\pi_{t'+1}(E^2(\G{Q(\beta)}))$ by the arguments given in Remark \ref{rem:category-W}.
This implies that $\iota^{\G{Q(\beta)}}_1{\comp}\epsilon_0$ is trivial in $\pi_{t'+1}(P^2(\G{Q(\beta)}))$.
This yields the following homotopy relative to $S^2$:
\begin{align*}&
P^2(\G{\hat{\gamma}_0}){\comp}\sigma_2(Q(\beta'))
\sim 
P^2(\G{\hat{\gamma}_0}){\comp}\sigma'_2(Q(\beta')) + P^2(\G{\hat{\gamma}_0}){\comp}\iota^{\G{Q(\beta')}}_1{\comp}\delta'_0
\\&\qquad
\sim 
P^2(\G{\hat{\gamma}_0}){\comp}\sigma'_2(Q(\beta')) + \iota^{\G{Q(\beta)}}_1{\comp}\delta_0{\comp}{\Sigma}\gamma + \iota^{\G{Q(\beta)}}_1{\comp}\epsilon_0
\\&\qquad
\sim 
\sigma'_2(Q(\beta)){\comp}\hat{\gamma}_0 + \iota^{\G{Q(\beta)}}_1{\comp}\delta_0{\comp}{\Sigma}\gamma
\sim 
\sigma_2(Q(\beta)){\comp}\hat{\gamma}_0.
\end{align*}
This completes the proof of Proposition \ref{prop:primitivity}.
\end{Proof*}

The attaching map $\psi(\beta)$ of the top cell of $E(\beta)$ is given by the composition
\begin{equation*}
\multiply\dgARROWLENGTH by2
\begin{diagram}
\node{S^t{\ast}S^1} \arrow{e,t}{[C(1_{S^{t}}),\Sigma{1_{S^1}}]^r} \node{C(S^{t}){\cup}S^t{\times}S^2} \arrow{e,t}{\hat{\chi}_{\beta}} \node{Q(\beta),}
\end{diagram}
\end{equation*}
where $\hat{\chi}_{\beta}$ denotes the map defined by $\hat{\mu}_{\beta}\vert_{S^t{\times}S^2} = \mu_0{\comp}(\beta{\times}1_{S^2})$ and $\hat{\mu}_{\beta}\vert_{C(S^t)} = \chi_{t+1}$ the characteristic map of the top cell of $Q(\beta)$.
Then a direct calculation shows that the following diagram is strictly commutative:
\begin{equation}\label{cd:h_2-psi}
\multiply\dgARROWLENGTH by2
\dgHORIZPAD=.5em
\dgVERTPAD=1ex
\begin{diagram}
\node{S^{t'}{\ast}S^1}
	\arrow{e,t}{[C(1_{S^{t'}}),\Sigma{1_{S^1}}]^r}
	\arrow{s,r}{\gamma{\ast}1_{S^1}}
\node{C(S^{t'}){\cup}S^{t'}{\times}S^2}
	\arrow{e,t}{\hat{\chi}_{\beta'}}
	\arrow{s,r}{\bar{\gamma}_0}
\node{Q(\beta')}
	\arrow{s,r}{\hat{\gamma}_0}
\\
\node{S^{t}{\ast}S^1}
	\arrow{e,t}{[C(1_{S^{t}}),\Sigma{1_{S^1}}]^r}
\node{C(S^t){\cup}S^t{\times}S^2}
	\arrow{e,t}{\hat{\chi}_{\beta}}
\node{Q(\beta)}
\end{diagram}
\end{equation}
where $\bar{\gamma}_0$ is given by $\bar{\gamma}_0\vert_{C(S^{t'})} = C(\gamma)$ and $\bar{\gamma}_0\vert_{S^{t'}{\times}S^2} = {\gamma}{\times}1_{S^2}$.
Thus we have that $\hat{\gamma}_0{\comp}\psi(\beta{\comp}\gamma) \sim \psi(\beta){\comp}(\gamma{\ast}1_{S^1})$.
By Proposition \ref{prop:fundamental-property} and Proposition \ref{prop:primitivity}, we have the following theorem.
\begin{Thm}\label{thm:unstable-eq}
For a map $\beta : S^{t} \to S^3$ and a suspension map $\gamma : S^{t'} \to S^{t}$ with $3 \leq t \leq t'$, we have that $E^3(\G{\hat{\gamma}_0})_{\ast}H^{SS}_2(\psi(\beta{\comp}\gamma)) = (\gamma{\ast}1_{S^1})^{\ast}H^{SS}_2(\psi(\beta)) = \pm({\Sigma}^2\gamma)^{\ast}H^{SS}_2(\psi(\beta))$
\end{Thm}
\begin{Cor}\label{cor:unstable-eq}
\begin{description}
\item{(1)} If $\beta{\comp}\gamma = 0$ in $\pi_{t'+1}(S^3)$, then $({\Sigma}^2\gamma)^{\ast}H^{SS}_2(\psi(\beta)) = \{0\}$.
\item{(2)} If $\beta : S^{t} \to S^3$ is of finite order $\ell$ with $t \geq 3$, then each element of $H^{SS}_2(\psi(\beta))$ is also of finite order which divides $\ell$.
\end{description}
\end{Cor}

We now prove the following lemma, making use of the notation of \cite{Toda:composition-methods}.

\begin{Lem}\label{thm:example}
Let $p$ be an odd prime, $\beta$ the co-H-map $\alpha_1(3) : S^{2p} \to S^3$ and $\gamma$ the suspension map $\alpha_2(2p) = \Sigma^{2p-3}\alpha_2(3) : S^{6p-5} \to S^{2p}$.
Then ${\Sigma}_{\ast}H^{S}_2(\psi(\alpha_1(3){\comp}\alpha_2(2p)))$ contains the composition of $\pm \alpha_1(6){\comp}\alpha_2(2p+3)$ with the bottom-cell inclusion.
\end{Lem}
\begin{Rem}
The composition $\alpha_1(5){\comp}\alpha_2(2p+2)$ is trivial for all odd primes $p$ except the prime $3$.
At the prime $3$, $\alpha_1(5){\comp}\alpha_2(8) = -3\beta_1(5) \not= 0$ and $\alpha_1(7){\comp}\alpha_2(10) = 0$ by Lemma 13.8 and Theorem 13.9 in \cite{Toda:composition-methods}.
\end{Rem}
%
\begin{Proof}
Firstly we summarise here some well-known results on odd primary components of stable and unstable homotopy groups of spheres.
\par
By Theorem 6.2 of Oka \cite{Oka:odd-component_1}, we know the following fact:
\begin{Fact}\label{fact:oka}
Let $p$ be an odd prime.
For $k < 2(p+3)(p-1)-4$, the $p$-component of stable homotopy group of the $k$-stem ${}_p\pi^S_k$ is trivial unless $k$ = $2r(p-1)-1$ ($1 \leq r \leq p+2$), $2p(p-1)-2$ and $2(p+1)(p-1)-3$.
In addition, all the non-trivial groups are given as follows: 
\begin{align*}&
{}_p\pi^S_{2r(p-1)-1} \cong \integer/p\integer\quad \text{generated by $\alpha_r$, $r \not= p$},
\\&
{}_p\pi^S_{2p(p-1)-1} \cong \integer/p^2\integer\quad \text{generated by $\alpha'_p$},
\\&
{}_p\pi^S_{2p(p-1)-2} \cong \integer/p\integer\quad \text{generated by $\beta_1$},
\\&
{}_p\pi^S_{2(p+1)(p-1)-3} \cong \integer/p\integer\quad \text{generated by $\alpha_1{\comp}\beta_1$},
\end{align*}
where $\alpha_r$ is defined inductively using Toda brackets: $\alpha_r = \langle\alpha_{r-1},p\iota,\alpha_1\rangle$.
Note that $\alpha_p = p\alpha'_p$.
At $p=3$, we remark that $\beta_1$ is given by a Toda bracket $\langle\alpha_{1},\alpha_{1},\alpha_{1}\rangle$.
\end{Fact}
\par
On the other hand, by (13.5) in \cite{Toda:composition-methods}, we know the following fact:
\begin{Fact}\label{fact:toda}\quad
$\pi_{2m-1+k}(S^{2m-1})_{(p)} \cong {}_p\pi^S_{k}$ if $\frac{k+3}{2(p-1)} \leq m$.
\end{Fact}
Since $\frac{2p(p-1)-2 + 3}{2(p-1)} < p+1$ and $2(p+1)-1 = 2p+1$, we have $\pi_{2(p+1)(p-1)+1}(S^{2p+1})_{(p)} \cong {}_p\pi^S_{2p(p-1)-2} \cong \integer/p\integer \ni \beta_1$.
Hence there is a generator $\beta_1(2p+1)$ of $\pi_{2(p+1)(p-1)+1}(S^{2p+1})_{(p)} \cong \integer/p\integer$ corresponding to the stable element $\beta_1$.
\par
In $\pi_{2(p+2)(p-1)+3}(S^{6})$, we know the following fact:
\begin{Fact}\label{eq:unstable-eq}\quad
$\alpha_2(6){\comp}\beta_1(4p+1) = 0$ in $\pi_{2(p+2)(p-1)+3}(S^6)_{(p)}$, 
which is obtained by a similar argument given in Page 184 of \cite{Toda:composition-methods} using (13.8) and Propositions 1.4 and 1.3 in \cite{Toda:composition-methods}:
\begin{align*}
\alpha_2(6){\comp}\beta_1(4p+1) & \in
2\{\alpha_1(6),\alpha_1(2p+3),p\iota_{4p}\}_1{\comp}\beta_1(4p+1)
\\&
=
2\alpha_1(6){\comp}\Sigma\{\alpha_1(2p+2),p\iota_{4p-1},\beta_1(4p-1)\}
\\&
\subset
2\alpha_1(6){\comp}\{\alpha_1(2p+3),p\iota_{4p},\beta_1(4p)\}_1,
\end{align*}
where $\{\alpha_1(2p+3),p\iota_{4p},\beta_1(4p)\}_1$ is a subset of $\pi_{2(p+2)(p-1)+3}(S^{2p+3})_{(p)} \cong {}_p\pi^S_{2(p+1)(p-1)-2} = 0$, since $\frac{2(p+1)(p-1)-2 + 3}{2(p-1)} < p+2$ and $2(p+2)-1 = 2(p-1)+5$.
\end{Fact}
\par
In $\pi_{6p-5}(S^{3})$, we also know the following fact:
\begin{Fact}\label{eq:unstable-eq3}\quad
$\alpha_1(3){\comp}\alpha_2(2p) = 2\alpha_2(3){\comp}\alpha_1(4p-2)$,
which is obtained by a similar argument given in Page 184 of \cite{Toda:composition-methods} using (3.9) (i) and Propositions 1.4 and 13.6 in \cite{Toda:composition-methods}:
\begin{align*}
\alpha_2(3){\comp}\alpha_1(4p-2) & \in \{\alpha_1(3),{\Sigma}(p\iota_{2p-1}),{\Sigma}\alpha_1(2p-1)\}_1{\comp}{\Sigma}^2\alpha_1(4p-4)
\\&\qquad
= \alpha_1(3){\comp}{\Sigma}\{p\iota_{2p-1},\alpha_1(2p-1),\alpha_1(4p-4)\},
\end{align*}
where the unstable Toda bracket 
$ 
\{p\iota_{2p-1},\alpha_1(2p-1),\alpha_1(4p-4)\} \subset \pi_{6p-6}(S^{2p-1})_{(p)} \cong {}_p\pi^S_{4p-5} \cong \integer/p\integer\{\alpha_{2}\}
$ 
corresponds to the stable Toda bracket $\langle{p\iota,\alpha_1,\alpha_1}\rangle = \frac{1}{2}\alpha_2$, since $\frac{4p-5+3}{2(p-1)} < 3 \leq p$ and $2p-1=2(p-1)+1$.
Thus $\{p\iota_{2p-1},\alpha_1(2p-1),\alpha_1(4p-4)\}$ determines $\frac{1}{2}\alpha_2(2p-1)$, 
and hence
$\alpha_1(3){\comp}\alpha_2(2p) = 2\alpha_2(3){\comp}\alpha_1(4p-2)$ in $\pi_{6p-5}(S^3)_{(p)}$.
\end{Fact}
Next we apply these facts to higher Hopf invariants.
\par
By Fact \ref{eq:unstable-eq} we have $\alpha_2(2p){\comp}\beta_1(6p-5) = 0$, since $p \geq 3$.
Then by Theorem \ref{thm:unstable-eq}, 
\begin{equation}\label{eq:unstable-eq2}
\begin{split}&
({\Sigma}^2\beta_1(6p-5))^{\ast}H^{SS}_2(\psi(\alpha_1(3){\comp}\alpha_2(2p)))
\\&\qquad
= \pm E^{3}(\G{\widehat{\beta_1(6p-5)}_0})_{\ast}H^{SS}_2(\psi(\alpha_1(3){\comp}\alpha_2(2p){\comp}\beta_1(6p-5))) = \{0\}.
\end{split}
\end{equation}
\par
By Fact \ref{eq:unstable-eq3} and Theorem \ref{thm:unstable-eq}, we have 
\begin{align*}&
E^{3}(\G{\widehat{\alpha_1(4p-2)}_0})_{\ast}H^{SS}_2(\psi(\alpha_1(3){\comp}\alpha_2(2p)))
\\&\qquad
= E^{3}(\G{\widehat{\alpha_1(4p-2)}_0})_{\ast}H^{SS}_2(\psi(2\alpha_2(3){\comp}\alpha_1(4p-2)))
\\&\qquad
= 2e({\Sigma}^2\alpha_1(4p-2))^{\ast}H^{SS}_2(\psi(\alpha_2(3))) 
= 2e\alpha_1(4p)^{\ast}H^{SS}_2(\psi(\alpha_2(3))), 
\end{align*}
where $e=\pm 1$.
Here we see that $E^{3}(\G{Q(\alpha_2(3))}) \simeq \G{Q(\alpha_2(3))}{\ast}\G{Q(\alpha_2(3))}{\ast}\G{Q(\alpha_2(3))}$ has the homotopy type of a wedge sum of spheres up to dimension $4p+1 (> 4p-2)$.
Since the suspensions of Whitehead products are trivial, each element of ${\Sigma}H^{SS}_2(\psi(\alpha_2(3)))$ has the form $a\iota_1{\comp}\alpha_2(6) + b\iota_2{\comp}\alpha_1(2p+4)$ by dimensional considerations, using Facts \ref{fact:oka}, \ref{fact:toda} and Corollary \ref{cor:unstable-eq} (2), where $a,b$ are integers modulo $p$ and $\iota_1 : S^{6} \hookrightarrow {\Sigma}E^{3}(\G{Q(\alpha_{2}(3))})$ and $\iota_2 : S^{2p+4} \hookrightarrow {\Sigma}E^{3}(\G{Q(\alpha_{2}(3))})$ denote appropriate inclusion maps.
Since $\alpha_1(5){\comp}\alpha_1(2p+2) = 0$ by Facts \ref{fact:oka} and \ref{fact:toda}, we have $\alpha_1(2p+4){\comp}\alpha_1(4p+1) = 0$, and hence each element of $2\alpha_1(4p+1)^{\ast}{\Sigma}_{\ast}H^{SS}_2(\psi(\alpha_2(3)))$ has the form $2a\alpha_2(6){\comp}\alpha_1(4p+1) = a\alpha_1(6){\comp}\alpha_2(2p+3)$ by Fact \ref{eq:unstable-eq3}.
On the other hand, we know
\begin{align*}&
2h_2(\eta_2{\comp}\alpha_2(3)){\comp}\alpha_1(4p-2) = 2h_2(\eta_2{\comp}\alpha_2(3){\comp}\alpha_1(4p-2)) 
\\&\qquad
= h_2(\eta_2{\comp}\alpha_1(3){\comp}\alpha_2(2p)) = \alpha_1(3){\comp}\alpha_2(2p) = 2\alpha_2(3){\comp}\alpha_1(4p-2),
\end{align*}
and hence ${\Sigma}^{3}(h_2(\eta_2{\comp}\alpha_2(3)){\comp}\alpha_1(4p-2)) = \alpha_2(6){\comp}\alpha_1(4p+1)$.
By the proof of Theorem \ref{thm:cat-E-2}, ${\Sigma}^{3}h_2(\eta_2{\comp}\alpha_2(3)){\comp}\alpha_1(4p-2)$ is given by the composition of $\pm\alpha_1(4p+1)^{\ast}{\Sigma}_{\ast}H^{SS}_2(\psi(\alpha_2(3)))$ with an appropriate inclusion map by Theorem \ref{thm:reduced-diagonal} and Remark \ref{rem:reduced-diagonal}.
Thus we have $a = \pm 1$ and 
$2\alpha_1(4p+1)^{\ast}{\comp}{\Sigma}_{\ast}H^{SS}_2(\psi(\alpha_2(3))) = \{\pm\alpha_1(6){\comp}\alpha_2(2p+3)\}$,
and hence we have 
\begin{equation}\label{eq:unstable-property}
\begin{split}&
{\Sigma}E^{3}(\G{\widehat{\alpha_1(4p-2)}_0})_{\ast}{\comp}{\Sigma}_{\ast}H^{SS}_2(\psi(\alpha_1(3){\comp}\alpha_2(2p))) 
\\&\qquad
= 2e\alpha_1(4p+1)^{\ast}{\comp}{\Sigma}_{\ast}H^{SS}_2(\psi(\alpha_2(3))) 
= \{e'\alpha_1(6){\comp}\alpha_2(2p+3)\},\quad e' = \pm 1.
\end{split}
\end{equation}
\par
By (13.10) in \cite{Toda:composition-methods}, we know that $h_p(\alpha_2(3)) = x\alpha_1(2p+1)$ for some $x \not= 0 \in \integer/p\integer$ and there are no other non-trivial James Hopf invariants $h_j(\alpha_2(3)), 1 < j \not= p$ for dimensional reasons.
Hence we have $\Sigma{\ad}(\alpha_2(3)) - \iota_1{\comp}\alpha_2(3) + \iota_1{\comp}\alpha_2(3) = H_1(\alpha_2(3)) = x\iota_2{\comp}\alpha_1(2p+1)$, where $\iota_1 : S^3 \to {\Sigma}\G{S^3}$ and $\iota_2 : S^{2p+1} \to {\Sigma}\G{S^3}$ denote appropriate inclusion maps (see Remark \ref{rems:(5)}).
Thus the attaching map of the $4p-1$ cell in $\Sigma\G{Q(\alpha_2(3))} - \Sigma\G(S^2)$ corresponding to that in $Q(\alpha_2(3))=S^2 \cup_{\eta_2{\comp}\alpha_2(3)} e^{4p-1}$ is given by
\begin{align}\label{align:H_1(alpha_2)}&
\Sigma{\ad}(\eta{\comp}\alpha_2(3))
= {\Sigma}\G{\eta}{\comp}\Sigma{\ad}(\alpha_2(3))
= {\Sigma}\G{\eta}{\comp}(\Sigma{\ad}(\alpha_2(3)) - \alpha_2(3) + \alpha_2(3))
\\&\qquad
= {\Sigma}\G{\eta}{\comp}(\Sigma{\ad}(\alpha_2(3)) - \alpha_2(3)) + {\Sigma}\G{\eta}{\comp}\alpha_2(3)
= x{\Sigma}\G{\eta}{\comp}\alpha_1(2p+1) + {\Sigma}\G{\eta}{\comp}\alpha_2(3).
\end{align}
Then it follows that the attaching maps of cells in 
$$
\G{Q(\alpha_2(3))}{\ast}\G{Q(\alpha_2(3))}{\ast}\G{Q(\alpha_2(3))} - \G(S^2){\ast}\G(S^2){\ast}\G(S^2)
$$
up to dimension $8p-2 (> 6p-2)$ are given by suspensions of ${\Sigma}\ad{(\eta{\comp}\alpha_2(3))}$.
Hence by (\ref{eq:unstable-property}) and (\ref{align:H_1(alpha_2)}), there is an integer $y$ such that each element of $\Sigma_{\ast}{H^{SS}_2(\psi(\alpha_1(3){\comp}\alpha_2(2p)))}$ has the form 
\begin{align*}&
e'\alpha_1(6){\comp}\alpha_2(2p+3) + y{\Sigma}^{2p}{\Sigma{\ad}(\eta{\comp}\alpha_2(3))}
\\&\qquad 
= e'\alpha_1(6){\comp}\alpha_2(2p+3) + y\{x{\Sigma}^{2p+1}\G{\eta}{\comp}\alpha_1(4p+1) + {\Sigma}^{2p+1}\G{\eta}{\comp}\alpha_2(2p+3)\}.
\end{align*}
If $y$ is non-zero modulo $p$, then each element of $\beta_1(6p-2)^{\ast}\Sigma_{\ast}{H^{SS}_2(\psi(\alpha_1(3){\comp}\alpha_2(2p)))}$ has the form 
\begin{align*}&
e'\alpha_1(6){\comp}\alpha_2(2p+3){\comp}\beta_1(6p-2) + y\{x{\Sigma}^{2p+1}\G{\eta}{\comp}\alpha_1(4p+1){\comp} + {\Sigma}^{2p+1}\G{\eta}{\comp}\alpha_2(2p+3)\}{\comp}\beta_1(6p-2)
\\&\qquad
= yx{\Sigma}^{2p}\G{\eta}{\comp}\alpha_1(4p+1){\comp}\beta_1(6p-2) \ne 0,
\end{align*}
since ${\Sigma}\G{S^2} \simeq {\Sigma}(S^1{\times}\G{S^3}) \simeq S^2 \vee {\Sigma}\G{S^3} \vee {\Sigma}^2\G{S^3}$.
This contradicts (\ref{eq:unstable-eq2}).
Thus we obtain $\Sigma_{\ast}{H^{SS}_2(\psi(\alpha_1(3){\comp}\alpha_2(4p-2)))} = \{e'\alpha_1(6){\comp}\alpha_2(2p+3)\}$, $e' = \pm 1$.
\end{Proof}
%

%
Using this, we show the following theorem.
\begin{Thm}\label{exam:punctured-m}
There is a series of simply connected closed $C^{\infty}$-manifolds $N_p$ indexed by odd primes $p \geq 5$ with $\cat{N_p} = \cat{(N_p - \{P\})}$, where $P$ is a point in $N_p$.
\end{Thm}
\begin{Rem}
The manifold $N_p$ does not have the property in Theorem \ref{exam:punctured-m} if we consider $\catq{}$, for any prime $q \ne p$, instead of $\cat{}$ or $\catp{}$.
\end{Rem}
\begin{Proof*}{\it Proof of Theorem \ref{exam:punctured-m}.}
We fix the prime $p \geq 5$ and let $L_p = E(\alpha_1(3){\comp}\alpha_2(2p))$ for the prime $p$ (see Theorem 13.4 in \cite{Toda:composition-methods}).
Then $L_p$ is a $C^{\infty}$-manifold with a CW decomposition $S^2 \cup_{\alpha} e^{6p-4} {\cup}_{\psi(\beta)} e^{6p-2}$, where $\alpha=\eta{\comp}\beta$ and $\beta = \alpha_1(3){\comp}\alpha_2(2p)$.
Here, $\alpha_1(3)$ is a co-H-map and $\alpha_2(2p)$ is a suspension map, and hence we have $h_2(\eta{\comp}\alpha_1(3){\comp}\alpha_2(2p)) = \alpha_1(3){\comp}\alpha_2(2p)$ by Proposition \ref{prop:fundamental-property}(1).
Also by Proposition 13.6 and (13.7) in \cite{Toda:composition-methods} and by the fact that $S^3$ is an H-space, we know that
\begin{equation*}
\alpha_1(3){\comp}\alpha_2(2p) \ne 0, \ \ {\Sigma}^1(\alpha_1(3){\comp}\alpha_2(2p)) \ne 0 \ \ \text{but} \ \ {\Sigma}^2(\alpha_1(3){\comp}\alpha_2(2p)) \in \pi_{6p-3}(S^{5})_{(p)} = 0,
\end{equation*}
which implies that $h_2(\eta{\comp}\alpha_1(3){\comp}\alpha_2(2p)) \ne 0$ and ${\Sigma}_{\ast}H^{S}_2(\psi(\alpha_1(3){\comp}\alpha_2(2p))) \ni 0$ by Lemma \ref{thm:example}.
Hence by Fact \ref{fact:cat-Q-1} and Theorem \ref{thm:cat-E-2}, we have $\cat{(S^2 {\cup}_{\alpha} e^{6p-4})} = 2$, $2 \leq \cat{L_p} \leq 3$.
If $\cat{L_p} = 2$, then we put $N_p=L_p$ which satisfies $\cat{(N_p - \{P\})} =  \cat{(S^2 {\cup}_{\alpha}e^{6p-4}) {\cup}_{\psi(\beta)} e^{6p-2}} = 2 = \cat{L_p}$.
Otherwise we have $\cat{L_p} = 3$ and then by Theorem \ref{thm:cat-E-3}, we also have $\cat{L_p{\times}S^n} = \cat {L_p} = 3$, $n\geq 1$, and then we put $N_p = L_p{\times}S^2$ which satisfies $\cat{(N_p - \{P\})} = \cat(L_p{\times}\{\ast\} \cup (S^2 {\cup}_{\alpha}e^{6p-4}){\times}S^2) = 3 = \cat{N_p}$.
Thus, in each case, there is a $C^{\infty}$-manifold which satisfies the required property.
\end{Proof*}

\begin{Thm}\label{exam:counter-ganea-m}
There is a simply connected closed $C^{\infty}$-manifold $M$ such that $\cat{M} = 3$ and $\cat{M{\times}S^n}$ $=$ $3$ for any $n \geq 2$ while we know only $3$ $\leq$ $\cat{M{\times}S^1}$ $\leq$ $4$ for $n=1$.
\end{Thm}
\begin{Cor}\label{exam2:counter-ganea-m}
There is a connected orientable closed $C^{\infty}$-manifold $N$ such that $\cat{N}$ $=$ $\cat{N{\times}S^n}$ for any $n \geq 1$.
\end{Cor}
\begin{Rems}
Ganea's conjecture for the manifold $M$ is true if we consider $\catp{}$, for any prime $p \ne 3$, instead of $\cat{}$ or $\catr{}$.
\end{Rems}

\begin{Proof*}{\it Proof of Theorem \ref{exam:counter-ganea-m}.}
Let $M = E(\alpha_1(3){\comp}\alpha_{2}(6))$ for the prime $p=3$ (see \cite{Toda:composition-methods}).
Then $M$ is a $C^{\infty}$-manifold with a CW decomposition $S^2 \cup_{\alpha} e^{14} {\cup}_{\psi(\beta)} e^{16}$, where $\alpha=\eta{\comp}\beta$ and $\beta = \alpha_1(3){\comp}\alpha_{2}(6)$.
Also by Theorem 13.4 in \cite{Toda:composition-methods} and by the fact that $S^5$ is an H-space at $p=3$, we know that
\begin{equation*}
\alpha_1(3){\comp}\alpha_{2}(6) \ne 0, \ \ {\Sigma}^3(\alpha_1(3){\comp}\alpha_{2}(6)) \ne 0 \in \pi_{16}(S^{6})_{(3)} \ \ \text{but} \ \ {\Sigma}^4(\alpha_1(3){\comp}\alpha_{2}(6)) \in \pi_{17}(S^{7})_{(3)} = 0,
\end{equation*}
which implies that ${\Sigma}^3_{\ast}h_2(\eta{\comp}\alpha_1(3){\comp}\alpha_2(6)) \ne 0$ and ${\Sigma}^2_{\ast}H^{S}_2(\psi(\alpha_1(3){\comp}\alpha_2(6))) \ni 0$ by Lemma \ref{thm:example}.
Hence by Theorems \ref{thm:cat-E-2} and \ref{thm:cat-E-3}, we have that $\cat{M} = 3$, $3 \leq \cat{M{\times}S^1} \leq 4$ and $\cat{M{\times}S^n} = 3$ for $n \geq 2$.
\end{Proof*}

\begin{Proof*}{\it Proof of Corollary \ref{exam2:counter-ganea-m}.}
Let $M = E(\alpha_1(3){\comp}\alpha_{2}(6))$ for the prime $p=3$ as in the proof of Theorem \ref{exam:counter-ganea-m}.
Then $M$ is a $C^{\infty}$-manifold with $\cat{M} = 3$, $3 \leq \cat{M{\times}S^1} \leq 4$ and $\cat{M{\times}S^n} = 3$ for $n \geq 2$.
If $\cat{M{\times}S^1} = 3$, we put $N=M$ which satisfies the required property.
Otherwise, we put $N=M{\times}S^{1}$ with $\psi'(\alpha_1(3){\comp}\alpha_{2}(6))$ as the attaching map of the top cell, where we know the set ${\Sigma}_{\ast}H^S_{3}(\psi'(\alpha_1(3){\comp}\alpha_{2}(6)))$ include the set ${\Sigma}^2_{\ast}H^S_{2}(\psi(\alpha_1(3){\comp}\alpha_{2}(6)))$ which contains ${\Sigma}^4\alpha_1(3){\comp}\alpha_{2}(6)=0$.
Then Theorem \ref{thm:obstruction-beta'} implies that $N$ satisfies the required properties.
\end{Proof*}
%
%

%
\end{document}